\documentclass[12pt]{amsart}
\usepackage{amssymb,amsmath,amsthm,mathrsfs}

\newtheorem{thm}{Theorem}
\newtheorem{lem}[thm]{Lemma}
\newtheorem{cor}[thm]{Corollary}
\newtheorem{prop}[thm]{Proposition}
\numberwithin{equation}{section}
\theoremstyle{definition}

\newtheorem{example}[thm]{Example}

\newtheorem*{remark}{Remark}

\newtheorem*{conj}{Conjecture}
\newcommand{\GL}{\mathrm{GL}}

\newcommand{\PGL}{\mathrm{PGL}}

\newcommand{\Res}{\mathrm{Res}}
\newcommand{\src}{\mathrm{src}}

\newcommand{\Ver}{\mathrm{Vert}}
\newcommand{\Edge}{\mathrm{Edge}}

\newcommand{\ga}{\alpha}
\newcommand{\gb}{\beta}

\newcommand{\Hom}{\mathrm{Hom}}

\newcommand{\Nul}{\mathrm{nullity}}
\newcommand{\Null}{\mathrm{Nul}}

\begin{document}
\title[Hecke Operators on Drinfeld Cusp Forms]{Hecke Operators on Drinfeld Cusp Forms}
\author{Wen-Ching Winnie Li and Yotsanan Meemark}
\address{Wen-Ching Winnie Li\\ Department of Mathematics\\ The Pennsylvania State University\\ University Park, PA 16802}
\email{\tt wli@math.psu.edu}
\address{Yotsanan Meemark\\ Department of Mathematics\\ Faculty of Science
\\ Chulalongkorn University\\ Bangkok, 10330 THAILAND}
\email{\tt yotsanan.m@chula.ac.th}
\thanks{The research of the first author is supported in part by the NSF grant DMS-0457574.
Part of the research was performed while she was visiting the
National Center for Theoretical Sciences, Mathematics Division, in
Hsinchu, Taiwan. She would like to thank the Center for its
support and hospitality. The second author was supported in part
by Grants for Development of New Faculty Staff from Chulalongkorn
University, Thailand.} 

\begin{abstract}
In this paper, we study the Drinfeld cusp forms for $\Gamma_1(T)$
and $\Gamma(T)$ using Teitelbaum's interpretation as harmonic
cocycles. We obtain explicit eigenvalues of Hecke operators
associated to degree one prime ideals acting on the cusp forms for
$\Gamma_1(T)$ of small weights and conclude that these Hecke
operators are simultaneously diagonalizable. We also show that the
Hecke operators are not diagonalizable in general for $\Gamma_1(T)$
of large weights, and not for $\Gamma(T)$ even of small weights. The
Hecke eigenvalues on cusp forms for $\Gamma(T)$ with small weights
are determined and the eigenspaces characterized.
\end{abstract}

\maketitle

\section{Introduction}

Hecke operators played a crucial role in the study of the
arithmetic of classical modular forms. Their actions on cusp forms
are skew Hermitian with respect to the Petersson inner product,
and hence they are diagonalizable. This property is fundamental in
understanding the classical cusp forms.

The function field analogue of the Poincare upper half plane is the
Drinfeld upper half plane. Parallel to the classical modular forms,
there are the Drinfeld modular forms introduced by Goss in
\cite{Gos80}. He also defined the Hecke operators in a similar way.
While certain arithmetic properties are alike for classical and
Drinfeld modular forms, there are also sharp differences. For
instance, B\"ockle \cite{Boc04} showed that the Eichler-Shimura
correspondence over a function field associates a Drinfeld
(cuspidal) common eigenform of Hecke operators to a degree one,
instead of degree two as in the classical case, Galois
representation, reflecting different multiplicative relations on
Hecke operators. Moreover, since the domain and image of Drinfeld
modular forms have the same positive characteristic, there is no
adequate analog of the Petersson inner product. Hence the
diagonalizability of the Hecke operators on Drinfeld forms still
remains an open question.

Using the residue map, Teitelbaum \cite{Tei91} in 1991 gave an
interpretation of Drinfeld cusp forms as harmonic cocycles on the
directed edges of a regular tree $\mathcal{T}$. The actions of the
Hecke operators were carried over to harmonic cocycles by B\"ockle
\cite{Boc04}. Since the directed edges of $\mathcal{T}$ are
parametrized by cosets of $\PGL_2$ over a local field $F$ modulo its
Iwahori subgroup $\mathcal I$, the Drinfeld cusp forms for a
congruence subgroup $\Gamma$ can then be regarded as vector-valued
left $\Gamma$-equivariant functions on $\PGL_2(F)/\mathcal I$, and
hence they are determined by the values on $\Gamma \backslash
\PGL_2(F)/\mathcal I$. This viewpoint is quite helpful in
computation when a fundamental domain is easily described. Another
advantage is that, by means of the strong approximation theorem, the
Drinfeld cusp forms can also be seen as equivariant functions in
adelic setting. This approach appeared in Gekeler and Reversat
\cite{GR96} and also in B\"ockle \cite{Boc04}.

Let $K =\mathbb{F}_q(T)$ be the rational function field. The
arithmetic of Drinfeld modular forms for the full modular group
$\GL_2(\mathbb{F}_q[T])$ was studied extensively in \cite{Gos80} and
\cite{Gek88}. Using geometric methods, B\"{o}ckle and Pink
investigated in \cite{Boc04} the structure of double cusp forms for
$\Gamma_1(T)$ with weight $k \le q+2$. They also computed the Hecke
eigenvalues for weight 4 double cusp forms.

The purpose of this paper is to study Drinfeld cusp and double cusp
forms for the congruence subgroups $\Gamma_1(T)$ and $\Gamma(T)$ of
$\GL_2(\mathbb{F}_q[T])$, with emphasis on the behavior of the Hecke
operators. Working with harmonic cocycles, we determine the
eigenvalues and the corresponding eigenspaces for Hecke operators at
degree one places of $K$. As we shall see, the diagonalizability of
the Hecke operators depends on the group and also the weight. More
precisely, the Hecke operators on the space of cusp forms of
$\Gamma_1(T)$ are diagonalizable for small weights $k \le q$, but
not for large weights $k > q$ in general. Further, as we pass from
$\Gamma_1(T)$ to its subgroup $\Gamma(T)$, the distinct eigenvalues
for Hecke operators on cusp forms with weights $k \le q$ remain the
same although the multiplicities may differ. We also characterize
each eigenspace. Explicit computations show that the Hecke actions
on the spaces of cusp forms and double cusp forms for $\Gamma(T)$ of
small weights change from diagonalizable to not diagonalizable as
the weight increases.

This paper is organized as follows. The Drinfeld cusp forms and
properties of the tree are reviewed in Sections 2 and 3,
respectively. Harmonic cocycles are recalled in Section 4. In
Section 5 we summarize Teitelbaum's isomorphism between Drinfeld
cusp forms and harmonic cocycles and describe B\"ockle's criterion
of double cusp forms as harmonic cocycles. The actions of the Hecke
operators on harmonic cocycles are introduced in Section 6. The body
of this paper is Sections 7 and 8, dealing with cusp forms for
$\Gamma_1(T)$ and $\Gamma(T)$, respectively. The final section gives
 examples of the Hecke actions on the cusp forms for
$\Gamma(T)$ for weights $k=3,4$ and 5, making explicit the main
results of the paper.

This paper grows out of the second author's thesis \cite{Mee06},
written under the direction of the first author.

The authors would like to thank the referee for valuable
comments and suggestions which corrected some errors and improved
the paper.

\medskip

\section{Drinfeld Cusp Forms}

Let $K=\mathbb F(T)$ be the rational function field over the finite
field $\mathbb{F}$ with $q$ elements. Write $\infty$ for the place
of $K$ with $1/T$ as a uniformizer. Then $A = \mathbb F[T]$ is the
ring of functions in $K$ regular outside $\infty$. Denote by
$K_\infty$ the completion of $K$ at $\infty$, $\mathcal{O}_\infty$
its ring of integers, and $\mathcal P_{\infty}$ the maximal ideal in
$\mathcal{O}_\infty$. Let $C = {{\widehat{\bar{K}}_\infty}}$ be the
completion of an algebraic closure of $K_\infty$.

The Drinfeld upper half plane $\Omega = C \smallsetminus K_\infty$
is endowed with a rigid analytic structure, on which
$\GL_2(K_\infty)$ acts by fractional linear transformations. For
$\gamma = \left(\begin{smallmatrix}
    a & b \\
    c & d \\
\end{smallmatrix}\right) \in \GL_2(K_\infty), m, k \in \mathbb{Z}$ and $f: \Omega \to
C$, define
\[
(f \underset{k,m}{|} \gamma)(z) := f(\gamma z) (\det \gamma)^m
(cz+d)^{-k}.
\]

Let $\Gamma$ be a congruence subgroup of the modular group
$\GL_2(A)$. It has finitely many cusps, represented by $\Gamma
\backslash \mathbb{P}^1(K)$. A rigid analytic function $f : \Omega
\to C$ is called a \emph{Drinfeld cusp form for $\Gamma$ of weight
$k$ and type $m$ for $\Gamma$} if it satisfies
\begin{itemize}
  \item[(i)] $f \underset{k,m}{|} \gamma = f$  for
  all $\gamma \in \Gamma$;
  \item[(ii)] $f$ is holomorphic at all cusps;
  \item[(iii)] $f$ vanishes at all cusps.
\end{itemize}
The cusp forms for $\Gamma$ of weight $k$ and type $m$ form a
vector space $S_{k,m}(\Gamma)$ over $C$. It contains a subspace
$S_{k,m}^2(\Gamma)$ of double cusp forms, which vanish at all
cusps at least twice.

\begin{remark}
While the weight can be any integer, the possible type is an
element in $\mathbb{Z}/(m_\Gamma)$, where $m_\Gamma$ is the order
of  $\det(\Gamma)$, a subgroup of $\mathbb{F}_q^\times$. Thus
$S_{k,m}(\Gamma) \ne 0$ implies $k \equiv 2m \mod (m_\Gamma)$. In
particular, if $m_\Gamma = 1$, which is the case to be considered
in this paper, then for fixed $k$, all $S_{k,m}(\Gamma)$ are
identical, and the same holds for $S_{k,m}^2(\Gamma)$.
\end{remark}

The following dimension formula for cusp forms was computed by
Teitelbaum.

\begin{prop}[\cite{Tei91}] \label{cdim}
Let $g_\Gamma$ be the genus of $\Gamma \backslash \bar{\Omega}$
and $h_\Gamma$ the number of cusps of $\Gamma \backslash \Omega$.
If $\Gamma$ is $p'$-torsion free and $m_\Gamma = 1$, then
\[
\dim_C S_{k,m} (\Gamma) = (k-1)(g_\Gamma + h_\Gamma - 1).
\]
\end{prop}

\medskip

\section{The Tree $\mathcal{T}$}

 The coset space
$\PGL_2(K_\infty)/\PGL_2(\mathcal{O}_\infty) =: \mathcal{T}$ may be
interpreted as a $(q+1)$-regular tree
 on which the group $\GL_2(K_\infty)$ acts by left
 translations. The vertices of $\mathcal{T}$ are the cosets
$\PGL_2(K_\infty)/\PGL_2(\mathcal{O}_\infty)$, while the directed
edges of $\mathcal{T}$ are parametrized by
$\PGL_2(K_\infty)/\mathfrak{I}_\infty$, where
$$\mathfrak{I}_\infty = \{ \left(\begin{smallmatrix}
    a & b \\
    c &  d \
     \end{smallmatrix}\right)
\in \GL_2(\mathcal{O}_\infty): c \in \mathcal P_{\infty} \}/\{
\left(\begin{smallmatrix}
    a & 0 \\
    0 &  a \
     \end{smallmatrix}\right)
\in \GL_2(\mathcal{O}_\infty)\}$$ is the Iwahori subgroup of
$\PGL_2(\mathcal{O}_\infty)$. The edge represented by $g \in
\GL_2(K_\infty)$ will be abbreviated as $\langle g \rangle$.

As in Serre \cite{Ser80}, a vertex or edge of $\mathcal{T}$ is
called \emph{$\Gamma$-stable} if its stabilizer in $\Gamma$ is
trivial; otherwise it is $\Gamma$-unstable. Let $\mathcal{T}_\infty$
be the subgraph of $\mathcal{T}$ consisting of unstable vertices and
edges. Then $S_0 = \Ver(\mathcal{T}) \smallsetminus
\Ver(\mathcal{T}_\infty)$ is the set of stable vertices and  $S_1 =
[\Edge(\mathcal{T}) \smallsetminus \Edge(\mathcal{T}_\infty)]/\pm$
is the set of non-oriented stable edges.

Two infinite paths in $\mathcal{T}$ are considered equivalent if they differ
at only finitely many edges. An $end$ of $\mathcal{T}$ is an
equivalence class of infinite paths $\{ e_1, e_2, \dots\}$. There is
a canonical bijection between the set of ends and
$\mathbb{P}^1(K_\infty)$, the boundary of $\Omega$; the rational
ends are $\mathbb{P}^1(K)$, corresponding to the cusps. The
stabilizer of an unstable vertex $v$ fixes a unique rational end,
and similarly for an unstable edge $e$; denote them by $b(v)$ and
$b(e)$, respectively. An edge $w$ of $\mathcal{T}$ is a
\emph{source} of an unstable edge $e$ if $w$ has the same
orientation as $e$ and there exists an unstable boundary vertex $v$
of $w$ such that the path from $v$ to its end $b(v)$ passes through
$e$. If $e$ is stable, then it is its own source. Denote by
$\src(e)$ the set of all sources of $e$. There are certain
inaccuracies in \cite{Tei91} concerning the sources of an edge. We
thank the referee for pointing them out.

\medskip

\section{Harmonic Cocycles}

For $k \ge 0$ and $m \in \mathbb{Z}$, let $V(k, m)$ be the
$(k-1)$-dimensional vector space over $C$ with a basis
$\{X^jY^{k-2-j}: 0 \le j \le k-2 \}$ endowed with the action of
$\GL_2(K_\infty)$ given by
\[ \gamma =
\left(\begin{smallmatrix}
    a & b \\
    c & d
\end{smallmatrix}\right) ~:
~X^jY^{k-2-j} \mapsto (\det \gamma)^{m-1}(dX-bY)^j(-cX+aY)^{k-2-j}
\]
for all $0 \le j \le k-2$. This then induces the action of $\gamma =
\left(\begin{smallmatrix}
    a & b \\
    c & d \
\end{smallmatrix}\right) \in \GL_2(K_\infty)$ on the dual space $\Hom(V(k,m),C)$
by sending $w \in \Hom(V(k,m),C)$ to
\[
(\gamma w )(X^jY^{k-2-j}) = (\det \gamma)^{1-m} w ((aX+bY)^j
(cX+dY)^{k-2-j})
\]
for $0 \le j \le k-2$.

A \emph{harmonic cocycle of weight $k$ and type $m$ for $\Gamma$} is
 a function $\mathbf{c}$ from the set of directed edges of $\mathcal{T}$
 to $\Hom(V(k, m), C)$ satisfying
        \begin{itemize}
            \item[(a)] For all vertices $v$ of $\mathcal{T}$, $$ \sum_{e \mapsto v}
\mathbf{c}(e) =
            0,$$ where $e$ runs through all edges in $\mathcal{T}$
            with terminal vertex $v$;
            \item[(b)] For all edges $e$ of $\mathcal{T}$, $\mathbf{c}(\bar{e})=-\mathbf{c}(e)$,
where $\bar{e}$ denotes $e$ with reversed orientation;
\item[(c)] It is $\Gamma$-equivariant, namely, for all
edges $e$ and elements $\gamma \in \Gamma$,
\[
\mathbf{c}(\gamma e) = \gamma(\mathbf{c}(e)).
\]
 \end{itemize}
The last condition means
\[
\mathbf{c}(\gamma e)(X^jY^{k-2-j}) = (\gamma
\mathbf{c}(e))(X^jY^{k-2-j}) = (\det \gamma)^{1-m} \mathbf{c}(e)
((aX+bY)^j (cX+dY)^{k-2-j})
\]
for all $\left(\begin{smallmatrix}
    a & b \\
    c & d
\end{smallmatrix}\right)\in \Gamma$ and $0 \le j \le k-2$.
Let $H_{k,m}(\Gamma)$ denote the space of harmonic cocycles of
weight $k$ and type $m$ for $\Gamma$.

As observed by Teitelbaum \cite{Tei91}, the value of a cocycle
$\mathbf{c} \in H_{k,m}(\Gamma)$ at a directed edge $e$ is the sum
of $\mathbf{c}$ evaluated at the source of $e$. Consequently,
cocycles in $H_{k,m}(\Gamma)$ are determined by their values on
$\Gamma \backslash S_1$.

\medskip

\section{Cusp Forms and Harmonic Cocycles}
 There is a building map from $\Omega$ to $\mathcal{T}$
commuting with the action of $\GL_2(K_\infty)$ (cf. \cite{Fv04} and
\cite{Tei91}). Using it one can define, for any $C$-valued
holomorphic 1-form $f(z)dz$ on $\Omega$, the residue $\Res_e f(z)dz$
at any directed edge $e$ of $\mathcal{T}$. This in turn gives a way
to associate harmonic cocycles to cusp forms. More precisely, for
each cusp
 form $f \in S_{k,m} (\Gamma)$, define the function $\Res(f)$ from the
 directed edges of $\mathcal{T}$ to $\Hom(V(k, m), C)$ by assigning,
 for any directed edge $e$, the values of $\Res(f)(e)$ at
the basis elements $X^jY^{k-2-j}$ to be
\begin{eqnarray}\label{Res}
\Res(f)(e)(X^jY^{k-2-j}) = \Res_e z^jf(z)dz
\end{eqnarray}
for all $0 \le j \le k-2$. Then properties (a) and (b) follow from
the rigid analytic residue theorem, and (c) from the modularity
 of $f$. Therefore $\Res(f)$ lies in $H_{k,m} (\Gamma)$.

\begin{thm} [Teitelbaum \cite{Tei91}] The residue map
$\mathrm{Res} : S_{k,m} (\Gamma) \to H_{k,m} (\Gamma)$ is an
isomorphism.
\end{thm}

Thus we identify cusp forms with harmonic cocycles. This allows us
to view cusp forms for $\Gamma$ as vector valued left
$\Gamma$-equivariant functions on
$\PGL_2(K_\infty)/\mathfrak{I}_\infty$, or left
$GL_2(K)$-equivariant functions on the adelic group $\GL_2(A_K)$
by applying the strong approximation theorem (cf. \cite{GR96} and
\cite{Rev00}). When $k = 2$, such functions are $C$-valued and
$\Gamma$-equivariance becomes $\Gamma$-invariance. Indeed, some
harmonic cocycles can be lifted to $\mathbb Z$-valued functions on
$\GL_2(K)\backslash \GL_2(\mathbb A_K)$, as remarked in
\cite{GR96},
 \cite{Rev00} and \cite{Boc04}.

Denote by $H_{k,m}^2(\Gamma)$ the image of $S_{k,m}^2(\Gamma)$ under
the $\Res$ map. To describe double cusp forms as cocycles, we define
the source of an end $[s]$ to be
\[
\src([s]) := \{ e : e \;\text{is stable, $t(e)$ is unstable and}\; b(t(e)) = [s] \},
\]
where $t(e)$ denotes the terminal vertex of $e$. The following
result of B\"ockle characterizes the image of double cusp forms
under the residue map.

\begin{thm} [B\"ockle \cite{Boc04}] \label{dchar}
Let $\Gamma_{[s]}$ denote the $\Gamma$-stabilizer of an end $[s]$
 representing a cusp of $\Gamma$. Then
\begin{itemize}
\item[(a)] The subspace of $V(k, m)$ stabilized by $\Gamma_{[s]}$,
denoted $V(k, m)^{\Gamma_{[s]}}$, is one-dimensional.
 \item[(b)] $\Gamma_{[s]}$ acts freely on
    $\src([s])$ with finitely many orbits, represented by edges
    $e_1^{[s]}$, $\dots, e_{l_s}^{[s]}$.
\item[(c)] Let $f \in S_{k,m}(\Gamma)$ and $\mathbf{c} = \Res(f)$.
 Then $f$ is a
double cusp form if and only if for any cusp $[s]$,
$\sum_{i=1}^{l_s} \mathbf{c}(e_i^{[s]})(g_s) = 0$ for any
generator $g_s$ of $V(k, m)^{\Gamma_{[s]}}$.
\end{itemize}
\end{thm}

Combined with Proposition \ref{cdim}, one obtains the dimension
formula for the space of double cusp forms:

\begin{prop} [B\"ockle \cite{Boc04}] \label{c2dim}
Let $g_\Gamma$ be the genus of $\Gamma \backslash \bar{\Omega}$
and $h_\Gamma$ the number of cusps of $\Gamma \backslash \Omega$.
If $\Gamma$ is $p'$-torsion free and $m_\Gamma = 1$, then
\[
\dim_C S_{k,m}^2 (\Gamma) = \left\{
                            \begin{array}{ll}
                              g_\Gamma & \hbox{if $k=2;$} \\
                              (k-2)(g_\Gamma + h_\Gamma -1) + g_\Gamma -1 & \hbox{if $k>2$.}
                            \end{array}
                          \right.
\]
\end{prop}

\medskip

\section{Hecke Operators}

We shall focus on the congruence groups $\Gamma = \Gamma_1(T)$ and
$\Gamma(T)$ defined as
\begin{equation*}
\Gamma_1(T) = \left\{ \left(\begin{smallmatrix}
    a & b \\
    c & d
\end{smallmatrix}\right) \in \GL_2(A) : a \equiv d \equiv 1 \; \text{and} \;
c \equiv 0 \mod T \right\}
\end{equation*}
and
\begin{equation*}
\Gamma(T) = \left\{ \left(\begin{smallmatrix}
    a & b \\
    c & d
\end{smallmatrix}\right) \in \GL_2(A) : a \equiv d \equiv 1 \; \text{and} \;
b \equiv c \equiv 0 \mod T \right\}.
\end{equation*}
They are $p'$-torsion free. Let $\mathfrak{P} \ne (T)$ be a maximal
ideal of $A$; choose the generator $P$ to be the irreducible
polynomial in $\mathfrak{P}$ satisfying $P(0) = 1$. Suppose $\deg P
= d$. Then
\[
\Gamma(T) \left(\begin{smallmatrix}
    P & 0 \\
    0 & 1
\end{smallmatrix}\right) \Gamma(T) = \Gamma(T)\left(\begin{smallmatrix}
    P & 0 \\
    0 & 1
\end{smallmatrix}\right) \sqcup \bigsqcup_{b \in A, \deg b < d} \Gamma(T)\left(\begin{smallmatrix}
    1 & b(1-P) \\
    0 & P
\end{smallmatrix}\right).
\]
The \emph{Hecke operator} at $\mathfrak{P}$ is defined using the
coset representatives of this double coset:
\[
T_\mathfrak{P} = P^{k-m-1} \bigl[ \left(\begin{smallmatrix}
    P & 0 \\
    0 & 1
\end{smallmatrix}\right) + \sum_{b \in A, \deg b < d} \left(\begin{smallmatrix}
    1 & b(1-P) \\
    0 & P
\end{smallmatrix}\right) \bigr],
\]
which acts on a holomorphic function $f$ on $\Omega$ via
$\underset{k,m}{|} T_\mathfrak{P}$. That is,
\begin{align*}
T_\mathfrak{P} f(z) = \bigl( f \underset{k,m}{|} T_\mathfrak{P}
\bigr)(z) = P^{k-m-1}\bigl[f
\underset{k,m}{|}\left(\begin{smallmatrix}
    P & 0 \\
    0 & 1
\end{smallmatrix}\right)(z) + \sum_{b \in A, \deg b < d} f
\underset{k,m}{|}\left(\begin{smallmatrix}
    1 & b(1-P) \\
    0 & P
\end{smallmatrix}\right)(z) \bigr].
\end{align*}

The generator $P$ is chosen in order to avoid the use of characters.
Here we have followed the normalization in B\"ockle \cite{Boc04},
which is a constant multiple of that defined by Goss \cite{Gos80}.
It is easy to check that $T_\mathfrak{P}$ sends $S_{k,m}(\Gamma)$ to
itself and preserves the
 double cusp forms. For two prime ideals $\mathfrak{P}$
and $\mathfrak{Q}$ not equal to $(T)$, $T_\mathfrak{P}$ commutes
with $T_\mathfrak{Q}$.

The action of the Hecke operator $T_\mathfrak{P}$ can be transported
to harmonic cocycles by means of the residue map. This was carried
out in \cite{Boc04}. Precisely, $T_\mathfrak{P}$ sends $\mathbf{c}
\in H_{k,m}(\Gamma)$ to a harmonic cocycle whose value at a directed
edge $e$ of $\mathcal{T}$ is
\begin{equation} \label{main}
\begin{aligned}
T_\mathfrak{P}\mathbf{c}(e) =
P^{k-m-1}\Bigl(\left(\begin{smallmatrix}
   P & 0 \\
   0 & 1
\end{smallmatrix}\right)^{-1}\!\mathbf{c}\left(\left(\begin{smallmatrix}
   P & 0 \\
   0 & 1
\end{smallmatrix}\right) e\right) + \sum_{b \in A, \deg b < d} \left(\begin{smallmatrix}
   1 & b(1-P) \\
   0 & P
\end{smallmatrix}\right)^{-1}\! \mathbf{c}\left(\left(\begin{smallmatrix}
   1 & b(1-P) \\
   0 & P
\end{smallmatrix}\right)e\right)\Bigr).
\end{aligned}
\end{equation}
This formula will be used to compute the eigenvalues and
eigenfunctions of Hecke operators. As we shall see from the cases
$\Gamma = \Gamma_1(T)$ and $\Gamma(T)$, the Hecke operators are
sometimes diagonalizable and sometimes not, depending on the group
and the weight.

\medskip

\section{Cusp Forms for $\Gamma_1(T)$}

In this section we consider cusp forms and double cusp forms for
$\Gamma_1(T)$. We may choose as a fundamental domain of $\Gamma_1(T)
\backslash \mathcal{T}$ the path connecting the cusp $[\infty] =
\tbinom{1}{0}$ and cusp $[0] = \tbinom{0}{1}$, as shown below. Recall from \S 3 
that $\langle g \rangle$ denotes the directed edge represented by $g$.
\[
[\infty] \quad \cdots \cdots \quad \leftrightarrows \left(\begin{smallmatrix}
    1 & 0 \\
    0 & T^2
\end{smallmatrix}\right) \underset{\left\langle\left(\begin{smallmatrix}
    1 & 0 \\
    0 & T^2
\end{smallmatrix}\right)\right\rangle}{\overset{\left\langle\left(\begin{smallmatrix}
    0 & 1 \\
    T & 0
\end{smallmatrix}\right)\right\rangle}{\genfrac{}{}{0pt}{1}{\longleftarrow}{\longrightarrow}}}
\left(\begin{smallmatrix}
    1 & 0 \\
    0 & T
\end{smallmatrix}\right) \underset{\bar{\gamma}_0 = \left\langle\left(\begin{smallmatrix}
    1 & 0 \\
    0 & T
\end{smallmatrix}\right)\right\rangle}{\overset{
\gamma_0 =\left\langle\left(\begin{smallmatrix}
    0 & 1 \\
    1 & 0
\end{smallmatrix}\right)\right\rangle}{
\genfrac{}{}{0pt}{1}{\longleftarrow}{\longrightarrow}}}
\left(\begin{smallmatrix}
    1 & 0 \\
    0 & 1
\end{smallmatrix}\right) \underset{\left\langle\left(\begin{smallmatrix}
    1 & 0 \\
    0 & 1
\end{smallmatrix}\right)\right\rangle}{\overset{\left\langle\left(\begin{smallmatrix}
    0 & T \\
    1 & 0
\end{smallmatrix}\right)\right\rangle}{\genfrac{}{}{0pt}{1}{\longleftarrow}{\longrightarrow}}}
\left(\begin{smallmatrix}
    T & 0 \\
    0 & 1
\end{smallmatrix}\right)  \leftrightarrows \quad \cdots \cdots \quad [0]
\]
It contains no stable vertices and one
stable edge $\langle\left(\begin{smallmatrix}
   0 & 1 \\
   1 & 0
\end{smallmatrix}\right)\rangle$, denoted by $\gamma_0$.
Then $g_{\Gamma_1(T)} = 0$ so that $\dim_CS_{k,m}(\Gamma_1(T))=
k-1$ by Proposition \ref{cdim}, and $\dim_C S_{2,m}^2(\Gamma_1(T))
= 0$ and $\dim_C S_{k,m}^2(\Gamma_1(T)) = k-3$ for $k \ge 3$ by
Proposition \ref{c2dim}. Theorem 3 of \cite{Tei91} implies
that any harmonic cocycle $\bf c$ for $\Gamma_1(T)$ automatically
vanishes on all edges of the fundamental domain except $\gamma_0$
and its two neighboring edges up to orientation. Further, the
value of $\bf c$ at $\gamma_0$ determines its values at the two
neighboring edges by harmonicity. Therefore to determine a
harmonic cocycle for $\Gamma_1(T)$, it suffices to first know its
value in Hom$(V(k,m), C)$ at $\gamma_0$, and then extend to other
edges by $\Gamma_1(T)$-equivariancy and harmonicity. This is
the strategy we shall use to compute the action of the Hecke
operators.

The stabilizers of the cusps $[\infty]$ and $[0]$ are
$(\Gamma_1(T))_{[\infty]} = \{\left(\begin{smallmatrix}
   1 & c \\
   0 & 1
\end{smallmatrix}\right) : c \equiv 0 \mod T \}$ and $(\Gamma_1(T))_
{[0]} = \{\left(\begin{smallmatrix}
   1 & 0 \\
   c & 1
\end{smallmatrix}\right) : c \equiv 0 \mod T \}$, respectively. Thus
$V(k,m)^{(\Gamma_1(T))_{[\infty]}}$ and
$V(k,m)^{(\Gamma_1(T))_{[0]}}$ are generated by $Y^{k-2}$ and
$X^{k-2}$, respectively. Also, $(\Gamma_1(T))_{[\infty]}
\backslash \src([\infty]) = \{\bar{\gamma}_0 = \langle
\left(\begin{smallmatrix}
   1 & 0 \\
   0 & T
\end{smallmatrix}\right)\rangle\}$ and $(\Gamma_1(T))_{[0]} \backslash \src([0]) = \{ \gamma_0
\}$. Recall that $\bar{\gamma}_0$ is the opposite of
$\gamma_0$. Hence by Theorem \ref{dchar}, we have

\begin{prop}
$S_{k,m}^2(\Gamma_1(T)) = \{ \mathbf{c} \in S_{k,m}(\Gamma_1(T)) :
\mathbf{c}(\gamma_0)(Y^{k-2}) =  \mathbf{c}(\gamma_0)(X^{k-2}) =
0\}.$
\end{prop}

Now we study the action of the Hecke operators $T_\mathfrak{P}$ on
$S_{k,m}(\Gamma_1(T))$, where $\mathfrak{P}$ is generated by $P =
1+\alpha T$. Using equation (\ref{main}), harmonicity and
$\Gamma_1(T)$-equivariancy, and noting $q$ is the cardinality of
the field $\mathbb F$, we get, for $0 \le j \le k-2$,
{\small\begin{equation}
\label{hecke}
\begin{aligned}
&T_\mathfrak{P} \mathbf{c}(\gamma_0) (X^jY^{k-2-j})= \mathbf{c}(\gamma_0)\Bigl( X^j(PY)^{k-2-j} \\
&+ \sum_{m=0}^{\lfloor\frac{j}{q-1}\rfloor} \Bigl(
\sum_{l=0}^{j-m(q-1)} \tbinom{j}{l+m(q-1)}\tbinom{k-2-j}{l}(1-P)^l -
P^{k-2-j} \tbinom{j}{m(q-1)} \Bigr) X^{j-m(q-1)} Y^{(k-2-j)+m(q-1)} \\
& + \sum_{n=1}^{\lfloor\frac{k-2-j}{q-1}\rfloor}
\Bigl(\sum_{l=n(q-1)}^{k-2-j} \tbinom{j}{l-n(q-1)}
\tbinom{k-2-j}{l}(1-P)^l - P^j \tbinom{k-2-j}{n(q-1)}
T^{n(q-1)}\Bigr) X^{j+n(q-1)} Y^{(k-2-j)-n(q-1)}\Bigr).
\end{aligned}
\end{equation}}

For each $0 \le j \le k-2$, define the harmonic cocycle
$\mathbf{c}_j$ by specifying its value at $\gamma_0$ by:
\begin{equation}\label{cj}
\mathbf{c}_j (\gamma_0) (X^jY^{k-2-j}) = 1 \quad {\rm and} \quad
~~~~\mathbf{c}_j (\gamma_0) (X^lY^{k-2-l}) = 0 \quad {\rm for} ~~l
\ne j. \end{equation} Further, put, for $0 \le j \le k-2$ and a
degree one polynomial $Q = 1+\beta T$, the polynomial
\begin{equation}\label{lambdaj}
\lambda_j(Q) = \sum_{l=0}^j \tbinom{j}{l} \tbinom{k-2-j}{l}(1-Q)^l =
\sum_{l=0}^{\min\{j, k-2-j\}} \tbinom{j}{l} \tbinom{k-2-j}{l}(-\beta
T)^l.
\end{equation}
Then $\lambda_j(Q)$ has degree at most ${\min\{j, k-2-j\}}$. Note
that $\lambda_0(Q) = \lambda_{k-2}(Q) = 1$ and $\lambda_j(Q) =
\lambda_{k-2-j}(Q)$ for all $0 \le j \le k-2$.

To see the behavior of the Hecke operators, we distinguish two
cases, according to the weight being small or large. First assume $q
\ge k \ge 2$. In this case (\ref{hecke}) is
 reduced to
\begin{equation} \label{qgek}
T_\mathfrak{P} \mathbf{c} (\gamma_0) (X^jY^{k-2-j}) = \lambda_j(P)
\mathbf{c} (\gamma_0) (X^j Y^{k-2-j}).
\end{equation}
Therefore each $\mathbf{c}_j$ is an eigenfunction of
$T_{\mathfrak{P}}$ with eigenvalue $\lambda_j(P)$. We have shown

\begin{thm} \label{cuspg1}
Let $\mathfrak{P}$ be a prime ideal of $A$ generated by
 $P$ with $P(0)=1$ and $\deg P = 1$. Suppose $q \ge k \ge 2$. Then
\begin{itemize}
    \item[$(1)$] Each $\mathbf{c}_j$, $0 \le j \le k-2$, is an
eigenfunction of the Hecke operator $T_\mathfrak{P}$ with
eigenvalue $\lambda_j(P)$; and
    \item[$(2)$] The Hecke operators at the ideals of degree one are simultaneously
diagonalized on $H_{k,m}(\Gamma_1(T))$ with respect to the basis
$\mathbf{c}_j$, $0 \le j \le k-2$.
\end{itemize}
\end{thm}

It is natural to ask if the $\mathbf{c}_j$, $0 \le j \le k-2$, are
also common eigenfunctions of the Hecke operators $T_\mathfrak{P}$
for prime ideals $\mathfrak{P}$ of degree $d >1$; and if so, find
the eigenvalues. Our computations lead to the following

\begin{conj} Let $\mathfrak{P}$ be a prime ideal of $A$ generated by
 $P$ with $P(0)=1$ and $\deg P = d \ge 1$. Suppose $q \ge k \ge 2$.
 Let $\theta$ be a root of $P$. Then each
$\mathbf{c}_j$, $0 \le j \le k-2$, is an eigenfunction of the Hecke
operator $T_\mathfrak{P}$ with eigenvalue $\lambda_j(P) :=
\prod_{i=0}^{d-1} \lambda_j(1 - \theta^{-q^i} T)$. Consequently, the
Hecke operators are simultaneously diagonalized on
$S_{k,m}(\Gamma_1(T))$.
\end{conj}

This conjecture is verified for $d \le 2$. Another evidence is
for the case $k=4$ and all $d$, provided by Prop. 15.6 in
\cite{Boc04}. It would be nice if the method there could be extended
to settle the conjecture.

\begin{remark}
If we factor the polynomial $\lambda_j (1+T) = \prod_{s=1}^{\deg
\lambda_j(1+T)} (1 + \delta_s T)$, then the eigenvalue $ \lambda_j
(P)$ above can also be expressed as $\prod_{s=1}^{\deg
\lambda_j(1+T)} P(\delta_s T)$.
\end{remark}

It is worth pointing out that the degree of $\lambda_j(P)$ above is
at most $d(k-2)/2$. This may be regarded as the Ramanujan conjecture
on Drinfeld cusp forms. A similar observation on weights can be
found in [Boc04], above Cor. 15.5.

\medskip

Notice that for $k \le q+2$ and $1 \le j \le q-2$, equation
(\ref{hecke}) is easily reduced to (\ref{qgek}) as well. Therefore
for $\mathfrak{P}$ of degree 1, $k=q+1$ and $1 \le j \le k-3$, one
gets
\[
T_\mathfrak{P} \mathbf{c}_j = \lambda_j(P) \mathbf{c}_j.
\]
Recall that a double cusp form $\mathbf{c}$ for $\Gamma_1(T)$
satisfies $\mathbf{c}(\gamma_0)(Y^{k-2}) =
\mathbf{c}(\gamma_0)(X^{k-2}) = 0$. Therefore $\mathbf{c}_1, \dots,
\mathbf{c}_{k-3}$ form a basis of the subspace of double cusp forms,
on which a similar result holds but with a slightly extended range
for $k$.

\begin{prop}
Let $\mathfrak{P}$ be a prime ideal of $A$ generated by the
polynomial $P$ of degree $1$ with $P(0)=1$.  If $q + 2 \ge k \ge
4$, then for all $\mathbf{c} \in S_{k,m}^2(\Gamma_1(T))$ and $1
\le j \le k-3$, one has
\[
T_\mathfrak{P} \mathbf{c} (\gamma_0) (X^jY^{k-2-j}) =
 \lambda_j(P)
\mathbf{c} (\gamma_0)( X^jY^{k-2-j}).
\]
\end{prop}
\begin{proof}
It remains to prove the proposition for the case $k = q+2$, and
$j=1$ or $k-3$. In this case, equation (\ref{hecke}) gives, for
$\mathbf{c} \in S_{k,m}^2(\Gamma_1(T))$, {\small\begin{align*}
T_\mathfrak{P} \mathbf{c}(\gamma_0)(XY^{k-3}) =\lambda_1(P)
\mathbf{c}(\gamma_0)(XY^{k-3})
\end{align*}}
and {\small\begin{align*} T_\mathfrak{P}
\mathbf{c}(\gamma_0)(X^{k-3}Y) =\lambda_{k-3}(P)
\mathbf{c}(\gamma_0)(X^{k-3}Y)
\end{align*}}
since $\mathbf{c}(\gamma_0)(Y^{k-2}) = \mathbf{c}(\gamma_0)(X^{k-2})
= 0$.
\end{proof}

\begin{cor} \label{dcuspg1}
Let $\mathfrak{P}$ be a degree one prime ideal of $A$ generated by
the polynomial $P$ with $P(0)=1$. If $q+2 \ge k \ge 4$, then
$\mathbf{c}_j$, $1 \le j \le k-3$, are eigenfunctions of the Hecke
operator $T_\mathfrak{P}$ with eigenvalue $\lambda_j(P)$. Further,
the Hecke operators for degree one prime ideals are simultaneously
diagonalized on $S_{k,m}^2(\Gamma_1(T))$ with respect to the basis
$\mathbf{c}_j$, $1 \le j \le k-3$.
\end{cor}

Note that there are no nonzero double cusp forms for weight $k < 4$.
The above result for $k = 4$ is Proposition 15.6 of \cite{Boc04},
proved by B\"{o}ckle and Pink.

\medskip

We now consider the case of general weight $k$. Assume $\mathfrak{P}
= (P)$, where $\deg P = 1$ and $P(0)=1$. Again, we appeal to
 (\ref{hecke}). For $i = 0, 1, \dots, q-2$ and $m_i = 0, 1,
\dots, \lfloor \frac{k-2-i}{q-1} \rfloor$, we have
\begin{align*}
&T_\mathfrak{P} \mathbf{c} (\gamma_0) (X^{i+m_i(q-1)}
Y^{k-2-(i+m_i(q-1))}) = \mathbf{c} (\gamma_0) \Bigl(X^{i+m_i(q-1)}(PY)^{k-2-(i+m_i(q-1))} \\
&+ \sum_{m=0}^{m_i} \Bigl( \sum_{l=0}^{i+(m_i-m)(q-1)}
\tbinom{i+m_i(q-1)}{l+m(q-1)}\tbinom{k-2-(i+m_i(q-1))}{l}(1-P)^l \\
&\quad \quad \quad \quad -P^{k-2-(i+m_i(q-1))}
\tbinom{i+m_i(q-1)}{m(q-1)} \Bigr)
X^{i+(m_i-m)(q-1)} Y^{k-2-(i+(m_i-m)(q-1))} \\
& + \sum_{n=1}^{\lfloor\frac{k-2-i}{q-1}\rfloor - m_i}
\Bigl(\sum_{l=n(q-1)}^{k-2-(i+m_i(q-1))}
\tbinom{i+m_i(q-1)}{l-n(q-1)}
\tbinom{k-2-(i+m_i(q-1))}{l}(1-P)^l \\
&\quad   - P^{i+m_i(q-1)} \tbinom{k-2-(i+m_i(q-1))}{n(q-1)}
T^{n(q-1)}\Bigr) X^{i+(m_i+n)(q-1)}
Y^{k-2-(i+(m_i+n)(q-1))}\Bigr).
\end{align*}
Recall the function $\mathbf{c}_j$ defined by (\ref{cj}). For $i =
0, 1, \cdots, q-2$, denote by $S_{k,m}(\Gamma_1(T))_i$ the subspace
of $S_{k,m}(\Gamma_1(T))$ generated by $\{ \mathbf{c}_i,
\mathbf{c}_{i+(q-1)}, \cdots, \mathbf{c}_{i+\lfloor
\frac{k-2-i}{q-1} \rfloor(q-1)} \}$ so that $S_{k,m}(\Gamma_1(T))
= \bigoplus_{i=0}^{q-2} S_{k,m}(\Gamma_1(T))_i$. The above calculation
proves the following
\begin{thm} \label{decom}
 Let $\mathfrak{P}= (P)$, where $\deg P = 1$ and $P(0) = 1$. Then for each $i = 0,$ $1,$ $\dots,$ $q-2$,
$S_{k,m}(\Gamma_1(T))_i$ is invariant under $T_\mathfrak{P}$. The
action of $T_\mathfrak{P}$ restricted to $S_{k,m}(\Gamma_1(T))_i$
with respect to the basis $\{ \mathbf{c}_i, \mathbf{c}_{i+(q-1)},
\cdots, \mathbf{c}_{i+\lfloor \frac{k-2-i}{q-1} \rfloor(q-1)} \}$ is
represented by the matrix $[T_\mathfrak{P}]_i = $
\[
{\tiny\left(\begin{matrix}
  \ga_{0,0}^{(i)}+P^{k-2-i} & \gb_{0,1}^{(i)} & \gb_{0,2}^{(i)} & \dots &
  \gb_{0,\lfloor\frac{k-2-i}{q-1}\rfloor}^{(i)} \\
  \ga_{1,1}^{(i)} & \ga_{1,0}^{(i)}+P^{k-2-(i+(q-1))} & \gb_{1,1}^{(i)} & \dots &
  \gb_{1,\lfloor\frac{k-2-i}{q-1}\rfloor - 1}^{(i)} \\
  \ga_{2,2}^{(i)} & \ga_{2,1}^{(i)} & \ga_{2,0}^{(i)}+P^{k-2-(i+2(q-1))} & \dots &
  \gb_{2,\lfloor\frac{k-2-i}{q-1}\rfloor - 2}^{(i)} \\
  \vdots & \vdots & \vdots & \ddots & \vdots \\
  \ga_{\lfloor\frac{k-2-i}{q-1}\rfloor,\lfloor\frac{k-2-i}{q-1}\rfloor}^{(i)}&
  \ga_{\lfloor\frac{k-2-i}{q-1}\rfloor,\lfloor\frac{k-2-i}{q-1}\rfloor-1}^{(i)}&
  \ga_{\lfloor\frac{k-2-i}{q-1}\rfloor,\lfloor\frac{k-2-i}{q-1}\rfloor-2}^{(i)}&
  \dots &
  \ga_{\lfloor\frac{k-2-i}{q-1}\rfloor,0}^{(i)}
  +P^{k-2-(i+\lfloor\frac{k-2-i}{q-1}\rfloor(q-1))}
\end{matrix}\right)}
\]
where {\small\begin{align*} \alpha_{m_i,m'}^{(i)} &=
\sum_{l=0}^{i+(m_i-m')(q-1)}
\tbinom{i+m_i(q-1)}{l+m'(q-1)}\tbinom{k-2-(i+m_i(q-1))}{l}(1-P)^l
-P^{k-2-(i+m_i(q-1))} \tbinom{i+m_i(q-1)}{m'(q-1)} \\
\intertext{and} \beta_{m_i,n}^{(i)} &=
\sum_{l=n(q-1)}^{k-2-(i+m_i(q-1))} \tbinom{i+m_i(q-1)}{l-n(q-1)}
\tbinom{k-2-(i+m_i(q-1))}{l}(1-P)^l \\
&-P^{i+m_i(q-1)} \tbinom{k-2-(i+m_i(q-1))}{n(q-1)} T^{n(q-1)}
\end{align*}}
for $m_i = 0, 1, \dots, \lfloor\frac{k-2-i}{q-1}\rfloor$, $0 \le m'
\le m_i$ and $1 \le n \le \lfloor\frac{k-2-i}{q-1}\rfloor - m_i$.
\end{thm}

Using geometric arguments, B\"{o}ckle and Pink computed the above
structures for the space of double cusp forms of $k=5, q=2$ and
$k=6,q=3$ in Proposition 15.3 of \cite{Boc04}. To illustrate the
above theorem, we give two examples of cusp forms with weights $k >
q$; in the first each Hecke action is diagonalizable, while in the
second it is not.

\begin{example} $q=3$, $k=7$ and $P=1+T$. There are two
invariant subspaces under $T_\mathfrak{P}$, namely,
$S_{7,m}(\Gamma_1(T))_0$ and $S_{7,m}(\Gamma_1(T))_1$ spanned by
$\{\mathbf{c}_0, \mathbf{c}_2, \mathbf{c}_4\}$ and
$\{\mathbf{c}_1, \mathbf{c}_3, \mathbf{c}_5\}$,
respectively. With respect to these bases, we have
\[
[T_\mathfrak{P}]_0 = \left(\begin{smallmatrix}
    1 & 0 &  0 \\
    2T^3 & 1 & T^3  \\
    2T & 2T & 1+2T
\end{smallmatrix}\right)
\; \mbox{and} \;  [T_\mathfrak{P}]_1 = \left(\begin{smallmatrix}
    1+2T & 2T^3 & 2 T^4 \\
    T & 1 & 2T^5 \\
    0 & 0 & 1
\end{smallmatrix}\right).
\]
Both matrices have the same distinct eigenvalues $1,  1+T+T\sqrt{1 -
T^2}$ and $1+T-T\sqrt{1 - T^2}$. Thus $[T_\mathfrak{P}]_0$ and
$[T_\mathfrak{P}]_1$ are diagonalizable, and hence so is
$T_\mathfrak{P}$.
\end{example}

\begin{example} \label{exgam1}
$q = 2$ and $k = 5$. There is only one polynomial $P=1+T$ to
consider. Further there is only one residue class mod $q-1$ given by
$i=0$, so one has
\[
[T_\mathfrak{P}]_0 = \left(\begin{smallmatrix}
   1   &  0 &  0   &  0  \\
  T^2  &  1 &  T^2 &  T^3 \\
  T    &  T &  1   &  T^3 \\
  0    &  0 &  0   &  1
\end{smallmatrix}\right).
\]
Thus $T_\mathfrak{P}$ has the eigenvalue $1$ of multiplicity two
with two linearly independent eigenfunctions $ \mathbf{c}_0$ and
$\mathbf{c}_3$, and the eigenvalue $1 + T^{3/2}$ of multiplicity two
with only one linearly independent eigenfunction
$T^{1/2}\mathbf{c}_1 + \mathbf{c}_2$. Hence $T_\mathfrak{P}$ is not
diagonalizable on $S_{5,m}(\Gamma_1(T))$. Further, since
$\mathbf{c}_1$ and $\mathbf{c}_2$ span the space of the double cusp
forms $S_{5,m}^2(\Gamma_1(T))$, this shows that the Hecke operator
$T_\mathfrak{P}$ is not diagonalizable on $S_{5,m}^2(\Gamma_1(T))$
either.
\end{example}

\begin{remark}
In both examples, unlike the case $k \le q+2$, there are
irrational eigenvalues. Our computations seem to suggest that the
nondiagonalizability results from inseparable eigenvalues. It would
be interesting to know if it could occur with separable eigenvalues.
\end{remark}

\medskip

\section{Cusp Forms for $\Gamma(T)$}

In this section, we work with $\Gamma = \Gamma(T)$, the group
of matrices in $\GL_2(A)$ congruent to the identity matrix modulo
$T$. A fundamental domain of $\Gamma(T) \backslash \mathcal{T}$ contains
$q+1$ rays, corresponding to the cusps $[\infty] = \tbinom{1}{0}$
and $[r] = \tbinom{r}{1}$, $r \in \mathbb{F}$, one stable vertex $
\left(\begin{smallmatrix}
   1 & 0 \\
   0 & 1
\end{smallmatrix}\right)$ and $q+1$ stable edges $ \gamma_r :=
\left\langle \left(\begin{smallmatrix}
   r & 1 \\
   1 & 0
\end{smallmatrix}\right)\right\rangle $, $ r \in \mathbb{F}$, and $\gamma_\infty :=
\left\langle \left(\begin{smallmatrix}
   1 & 0 \\
   0 & 1
\end{smallmatrix}\right)\right\rangle$. Thus $g_{\Gamma(T)} = 0$ so
that $\dim_C S_{k,m}(\Gamma(T)) = (k-1)q$ by Proposition \ref{cdim},
and $\dim_C S_{2,m}^2(\Gamma(T)) = 0$ and $\dim_C
S_{k,m}^2(\Gamma(T)) = (k-2)q - 1$ for $k \ge 3$ by Proposition
\ref{c2dim}. To determine a harmonic cocycle for $\Gamma(T)$, as
noted in Sec. 4, one needs to know only its values at $\gamma_r$, $r
\in \mathbb{F}$, and its value at $\gamma_\infty$ is determined by
the hamornicity condition $\mathbf{c}(\gamma_\infty) + \sum_{r \in
\mathbb{F}} \mathbf{c}(\gamma_r) = 0$. The stabilizer of the cusp
$[\infty]$ (resp. $[r]$, $r \in \mathbb{F}$) is $\Gamma_{[\infty]} =
\left\{ \left(\begin{smallmatrix}
   1 & c \\
   0 & 1
\end{smallmatrix}\right) : c \equiv 0 \mod T \right\}$ (resp. $\Gamma_{[r]} = \left\{
\left(\begin{smallmatrix}
   1+rc & -r^2c \\
   c    & 1-rc
\end{smallmatrix}\right) : c \equiv 0 \mod T \right\}$) so that $V(k,m)^{\Gamma_{[\infty]}}$
(resp. $V(k,m)^{\Gamma_{[r]}}$) is spanned by $ Y^{k-2} $ (resp.
$(X-rY)^{k-2}$). Moreover, $\Gamma_{[\infty]} \backslash
\src([\infty]) = \{ \gamma_\infty \}$ and $\Gamma_{[r]} \backslash
\src([r]) = \{ \gamma_r \}$, $r \in \mathbb{F}$. Thus by Theorem
\ref{dchar}, the double cusp forms can be described as follows.

\begin{prop} A harmonic cocycle $\mathbf{c} \in
H_{k,m}(\Gamma(T))$ lies in $H_{k,m}^2(\Gamma(T))$ if and only if $
\mathbf{c}(\gamma_\infty)(Y^{k-2}) = 0$ and
$\mathbf{c}(\gamma_r)((X- rY)^{k-2}) = 0$ for all $ r \in
\mathbb{F}.$
\end{prop}

Next we study the action of the Hecke operator $T_\mathfrak{P}$ at
$\mathbf{c} \in H_{k,m}(\Gamma(T))$. Recall that a harmonic cocycle
takes values in $\Hom(V(k,m), C)$. In view of the above proposition,
it turns out that the action is best described if, for all $r \in
\mathbb F$, the basis $(X-rY)^jY^{k-2-j}$, $0 \le j \le k-2$, of
$V(k,m)$ is used when we discuss the values of a harmonic cocycle at
the directed edge $\gamma_r$. Therefore we shall describe the action
using such bases. To ease our notation, for $\mathbf{c} \in
H_{k,m}(\Gamma(T))$, $r \in \mathbb F$, and $0 \le j \le k-2$, let
\begin{equation}\label{Z}
Z(\mathbf{c}, r, j) = \mathbf{c}(\gamma_r)((X-rY)^jY^{k-2-j}).
\end{equation}
Assume that $\mathfrak{P}$ is generated by $P = 1 + \alpha T$ with
$\alpha \in \mathbb F^\times$.  Again, we use (\ref{main}){,
harmonicity and $\Gamma(T)$-equivariancy to arrive at the main
identity of the Hecke action:}
\begin{equation}  \label{heckeg}
\begin{aligned}
Z(T_\mathfrak{P}\mathbf{c}, r, j)& = P^{k-2-j}Z(\mathbf{c}, r, j) -
P^j \sum_{n=1}^{\lfloor \frac{k-2-j}{q-1}
\rfloor} \tbinom{k-2-j}{n(q-1)} T^{n(q-1)} Z(\mathbf{c}, r,j+ n(q-1))  \\
 + \sum_{b \ne r} & \Bigl[ \sum_{u=0}^j (b-r)^{j-u} \Bigl(
P^{k-2-j} \tbinom{j}{u} - \sum_{l=0}^{u}
\tbinom{j}{u-l} \tbinom{k-2-j}{l} (1-P)^l \Bigr) Z(\mathbf{c},b,u) \\
& - \sum_{u=j+1}^{k-2} \sum_{l=u-j}^{k-2-j} \tbinom{j}{u-l}
\tbinom{k-2-j}{l} (1-P)^l (b-r)^{j-u} Z(\mathbf{c},b,u) \Bigr].
\end{aligned}
\end{equation}
Notice that when $j=k-2$, (\ref{heckeg}) becomes
\begin{equation} \label{gamk-2}
Z(T_\mathfrak{P}\mathbf{c},r, k-2) = Z(\mathbf{c},r,k-2)
\end{equation}
for all $r \in \mathbb{F}$. Moreover, for $j = 0$ and $r \in
\mathbb{F}$ we have
\begin{align*}
Z(T_\mathfrak{P} \mathbf{c}, r, 0) = & \: P^{k-2} Z(\mathbf{c},r,0)
- \sum_{n=1}^{\lfloor \frac{k-2}{q-1} \rfloor}
\tbinom{k-2}{n(q-1)} T^{n(q-1)} Z(\mathbf{c},r, n(q-1)) \\
& + \sum_{b \ne r} \Bigl( (P^{k-2} - 1) Z(\mathbf{c},b,0) -
\sum_{u=1}^{k-2} (1-P)^u (b-r)^{-u} Z(\mathbf{c},b,u) \Bigr).
\end{align*}
Summing over all $r \in \mathbb{F}$ and using harmonicity, we get
\begin{align*}
-&T_\mathfrak{P}  \mathbf{c} (\gamma_\infty) (Y^{k-2})  = I + II,
\end{align*}
where $$ I= \sum_{r \in \mathbb{F}} \Bigl( P^{k-2} Z(\mathbf{c},r,0)
- \sum_{n=1}^{\lfloor \frac{k-2}{q-1} \rfloor} \tbinom{k-2}{n(q-1)}
T^{n(q-1)} Z(\mathbf{c},r, n(q-1)) \Bigr)$$ and
\begin{eqnarray*}
II &=  \sum_{b \in \mathbb{F}}\Bigl( \sum_{r \ne b} (P^{k-2} - 1)
Z(\mathbf{c},b,0) - \sum_{u=1}^{k-2} (1-P)^u \sum_{r \ne b}
(b-r)^{-u} Z(\mathbf{c},b,u) \Bigr) \\
&= \sum_{b \in \mathbb{F}}\Bigl( - (P^{k-2} - 1) Z(\mathbf{c},b,0) +
\sum_{n=1}^{\lfloor \frac{k-2}{q-1}\rfloor} (1-P)^{n(q-1)}
Z(\mathbf{c},b,n(q-1)) \Bigr).
\end{eqnarray*}
Combined, this gives
\begin{equation} \label{gamk-2i}
T_\mathfrak{P} \mathbf{c} (\gamma_\infty) (Y^{k-2}) = \mathbf{c}
(\gamma_\infty) (Y^{k-2}).
\end{equation}
The equations (\ref{gamk-2}) and (\ref{gamk-2i}) then imply

\begin{prop} \label{eigenone} Let $\mathbf{c} \in S_{k,m}(\Gamma(T))$ be an eigenfunction
of $T_\mathfrak{P}$, where $\mathfrak{P} \ne (T)$ has degree 1. If
it is not a double cusp form, then the eigenvalue is $1$.
\end{prop}

Assume further that $q \ge k \ge 2$. In this case (\ref{heckeg})
is reduced to
\begin{equation} \label{redhegam}
\begin{aligned}
Z(T_\mathfrak{P} \mathbf{c}, r, j)
 &=
\sum_{u=0}^{j-1} \sum_{b \in \mathbb{F}} (b-r)^{j-u} \Bigl(
P^{k-2-j} \tbinom{j}{u} - \sum_{l=0}^{u} \tbinom{j}{u-l}
\tbinom{k-2-j}{l} (1-P)^l \Bigr) Z(\mathbf{c},b,u) \\
& \quad+  [P^{k-2-j} - \lambda_j(P)] \sum_{b \in \mathbb{F}}
Z(\mathbf{c},b,j) + \lambda_j(P) Z(\mathbf{c},r,j)\\
& \quad - \sum_{u=j+1}^{k-2} \sum_{l=u-j}^{k-2-j} \tbinom{j}{u-l}
\tbinom{k-2-j}{l} (1-P)^l \sum_{b \ne r}  (b-r)^{j-u} Z(\mathbf{c},b,u)\\
&= \sum_{u=0}^{j-1} \alpha_u(j, P) \sum_{b \in \mathbb{F}}
(b-r)^{j-u} Z(\mathbf{c},b,u)  +  [P^{k-2-j} - \lambda_j(P)]
\sum_{b \in \mathbb{F}} Z(\mathbf{c},b,j) \\
&+ \lambda_j(P) Z(\mathbf{c},r,j) - \sum_{u=j+1}^{k-2} \beta_u(j, P)
\sum_{b \ne r} (b-r)^{j-u} Z(\mathbf{c},b,u),
\end{aligned}
\end{equation}
where $\lambda_j(P) = \sum_{l=0}^j \tbinom{j}{l}
\tbinom{k-2-j}{l}(1-P)^l$ is given by (\ref{lambdaj}),
$$\alpha_u(j,P) =
P^{k-2-j}\tbinom{j}{u} - \sum_{l=0}^u \tbinom{j}{u-l}
\tbinom{k-2-j}{l}(1-P)^l \quad \text{for} \quad 0 \le u \le j-1,$$
and
$$\beta_u(j,P) = \sum_{l=u-j}^{k-2-j} \tbinom{j}{u-l}
\tbinom{k-2-j}{l}(1-P)^l \quad \text{for} \quad j+1 \le u \le
k-2.$$
 For $r
\in \mathbb{F}$ and $0 \le j \le k-2$, denote by
$\mathbf{c}_j^{(r)}$ the function
$$\mathbf{c}_j^{(r)} (\gamma_r)
((X-rY)^j Y^{k-2-j}) = 1 \quad \text{and} \quad \mathbf{c}_j^{(r)}
(\gamma_s) ((X-rY)^l Y^{k-2-l}) = 0 ~~\text{ if} ~~s \ne r
~~\text{ or} ~~l \ne j.$$ Let $\mathbf{c}_j = \sum_{r \in
\mathbb{F}} \mathbf{c}_j^{(r)}$. Then $T_\mathfrak{P}\mathbf{c}_j
= \lambda_j(P) \mathbf{c}_j$, that is, $\mathbf{c}_j$ is an
eigenfunction of $T_\mathfrak{P}$ with eigenvalue $\lambda_j(P)$.
Observe that $\mathbf{c}_j$ are liftings of the eigenfunctions of
$S_{k,m}(\Gamma_1(T))$.

Our next goal is to show that $\lambda_j(P)$ are the eigenvalues for
the Hecke operator $T_\mathfrak{P}$ on $S_{k,m}(\Gamma(T))$ when $q
\ge k$. For this, we need
\begin{lem} \label{red}
Suppose that $\mathbf{c}$ is an eigenfunction of the Hecke
operator $T_\mathfrak{P}$ on $S_{k,m}(\Gamma(T))$ with eigenvalue
$\lambda \ne \lambda_n(P)$ for all $0 \le n \le k-2$. Then for
each $0 \le n \le k-2$ and each $r \in \mathbb{F}$, there are
constants $A_u^{(n)} \in \mathbb{F}(T)$ for $n+1 \le u \le k-2$
such that
\begin{subequations} \label{mainlem}
\begin{equation} \tag*{(\ref{mainlem})$_n$}
(\lambda - \lambda_n(P)) Z(\mathbf{c}, r, n)
 = \sum_{u=n+1}^{k-2} A_u^{(n)} \sum_{b \ne r} (b-r)^{n-u}
Z(\mathbf{c}, b, u).
\end{equation}
\end{subequations}
\end{lem}

Grant this lemma. By applying (\ref{mainlem})$_n$ repeatedly from $n
= k-2$ down to $n=0$, we deduce that $\mathbf{c} = 0$. This proves

\begin{thm}\label{eigen}
Let $\mathfrak{P} =(P) \ne (T)$ be a degree one prime ideal of $A$.
For $q \ge k \ge 2$ the distinct eigenvalues for the Hecke operator
$T_\mathfrak{P}$ on $S_{k,m}(\Gamma(T))$ are the distinct
$\lambda_j(P)$, $0 \le j \le k-2$.
\end{thm}

Let $\mathbf{c}$ be an eigenfunction of $T_\mathfrak{P}$ with
eigenvalue $\lambda$. Then (\ref{redhegam}) gives rise to
{\small\begin{align*} (\lambda - \lambda_j(P)) Z(\mathbf{c},r,j)
=\: &\sum_{u=0}^{j-1} \alpha_u(j, P)\sum_{b \in \mathbb
F}(b-r)^{j-u} Z(\mathbf{c},b,u)  + \sum_{b \in \mathbb{F}}
[P^{k-2-j} - \lambda_j(P)]
Z(\mathbf{c},b,j) \\
&  -  \sum_{u=j+1}^{k-2} \beta_u(j, P) \sum_{b \ne r} (b-r)^{j-u}
Z(\mathbf{c},b,u)
\end{align*}}
for all $0 \le j \le k-2$ and $r \in \mathbb{F}$. Summing over all
$r \in \mathbb{F}$, we get, for each $0 \le j \le k-2$,
\begin{equation}\label{sumzero}
(\lambda - \lambda_j(P)) \sum_{r \in \mathbb{F}} Z(\mathbf{c},r,j) =
0.
\end{equation}
Hence if $\lambda \ne \lambda_j(P)$, then $\sum_{r \in \mathbb{F}}
Z(\mathbf{c},r,j) = 0$ so that
\begin{subequations} \label{pflemma}
\begin{equation}\tag*{(\ref{pflemma})$_j$}
\begin{aligned}
(\lambda - \lambda_j(P)) Z(\mathbf{c},r,j)  =\: &\sum_{u=0}^{j-1}
\alpha_u(j, P)\sum_{b \in \mathbb F}(b-r)^{j-u} Z(\mathbf{c},b,u)
\\
&- \sum_{u=j+1}^{k-2} \beta_u(j, P) \sum_{b \ne r} (b-r)^{j-u}
Z(\mathbf{c},b,u)
\end{aligned}
\end{equation}
\end{subequations}
for all $0 \le j \le k-2$ and $r \in \mathbb{F}$. When $j=0$, the
first sum on the right side is void and hence (\ref{mainlem})$_0$
holds with $A_u^{(0)} = \beta_u(0,P)$ for $1 \le u \le k-2$.  We
shall prove Lemma \ref{red} by induction on $n$. To proceed, we
prove an identity which will be used repeatedly in the
computations to follow.

\begin{prop} \label{intersum}
For $1 \le l,t \le k-2 \le q-2$ and any $C$-valued function $X(s)$
 on $\mathbb F$, we have
$$\sum_{b \in \mathbb{F}}\sum_{\substack{s \in \mathbb{F} \\ s \ne
b}} \dfrac{(b-r)^t}{(s-b)^l}X(s) = \left\{\begin{array}{ll}
   \sum_{s \in \mathbb{F}} (-1)^{l+1}\tbinom{t}{l}(s-r)^{t-l}X(s) & \hbox{$\mathrm{if}$ $t > l;$} \\
   \sum_{s \in \mathbb{F}} (-1)^{l+1}X(s) & \hbox{$\mathrm{if}$ $t = l;$} \\
   0 & \hbox{$\mathrm{if}$ $t \le l$}.
          \end{array}
   \right.$$
\end{prop}

\begin{proof}
Let $1 \le l,t \le k-2 \le q-2$. Then
\begin{align*}
\sum_{b \in \mathbb{F}}\sum_{\substack{s \in \mathbb{F} \\
s \ne b}} \frac{(b-r)^t}{(s-b)^l}X(s) &= \sum_{s \in
\mathbb{F}}\sum_{\substack{b \in \mathbb{F} \\ b \ne s}}
\frac{((b-s) + (s-r))^t}{(s-b)^l}X(s) \\  &=  \sum_{s \in
\mathbb{F}} \sum_{i=0}^t  (-1)^l \tbinom{t}{i} (s-r)^{t-i} \sum_{b
\ne s} (b-s)^{i-l}X(s).
\end{align*}
Since $1 \le l,t \le k-2 \le q-2$, $\sum_{\substack{b \in \mathbb
F \\ b \ne s}} (b-s)^{i-l}$ vanishes unless $i=l$ in which case it
is $-1$, so
\begin{align*}
\sum_{b \in \mathbb{F}}\sum_{\substack{s \in \mathbb{F} \\
s \ne b}} \frac{(b-r)^t}{(s-b)^l}X(s) =  &\left\{\begin{array}{ll}
   \sum_{s \in \mathbb{F}} (-1)^{l+1}\tbinom{t}{l}(s-r)^{t-l}X(s), & \hbox{$\mathrm{if}$ $t > l;$} \\
   0, & \hbox{$\mathrm{if}$ $t < l$}.
          \end{array}
   \right.
\end{align*}
If $t=l$, then
\begin{align*}
\sum_{b \in \mathbb{F}}\sum_{\substack{s \in \mathbb{F} \\
s \ne b}} &\Bigr(\frac{b-r}{s-b}\Bigr)^t X(s)= \sum_{s \in
\mathbb{F}}\sum_{\substack{b \in \mathbb{F} \\ b \ne s}}
\Bigr(\frac{s-r}{s-b} - 1\Bigr)^tX(s) \\
&= \sum_{s \in \mathbb{F}}\sum_{b \ne s} \Bigl( \sum_{i=0}^{t-1}
\tbinom{t}{i} (-1)^i \Bigr(\frac{s-r}{s-b}\Bigr)^{t-i} + (-1)^t
\Bigr)X(s) = (-1)^{t+1}\sum_{s \in \mathbb F} X(s).
\end{align*}
This proves the proposition.
\end{proof}

\begin{proof}[Proof of Lemma \ref{red}]
We shall apply Proposition \ref{intersum} to $X(s) = Z(\mathbf{c},s,
j)$, in which case the sum is equal to 0 when $t = l$ because of
(\ref{sumzero}) and the assumption $\lambda \ne \lambda_j(P)$ for
all $j$. Assume that the statement is valid up to $n$, where $0 \le
n < k-2$. That is, for all $0 \le j \le n$ and $b \in \mathbb{F}$,
we have
\begin{subequations} \label{indhypo}
\begin{equation} \tag*{(\ref{indhypo})$_j$}
Z(\mathbf{c},b,j) = \frac{1}{\lambda -
\lambda_j(P)}\sum_{u=j+1}^{k-2} A_u^{(j)} \sum_{s \ne b} (s-b)^{j-u}
Z(\mathbf{c},s,u).
\end{equation}
\end{subequations}
Substituting (\ref{indhypo})$_0$ into (\ref{pflemma})$_{n+1}$, we
get
{\small\begin{align*}
&(\lambda - \lambda_{n+1}(P)) Z(\mathbf{c},r, n+1) \\
& = \sum_{u=0}^{n} \alpha_u(n+1,P)  \sum_{b \in \mathbb{F}}
(b-r)^{(n+1)-u}  Z(\mathbf{c},b,u) - \sum_{u=n+2}^{k-2}
\beta_u(n+1,P) \sum_{b
\ne r}(b-r)^{(n+1)-u} Z(\mathbf{c},b,u) \\
& =  \sum_{b \in \mathbb{F}} \tfrac{\alpha_0(n+1,P)}{\lambda -
\lambda_0(P)} \sum_{u=1}^{k-2} A_u^{(0)} \sum_{s \ne b}
\tfrac{(b-r)^{n+1}}{(s-b)^u} Z(\mathbf{c},s,u)  + \sum_{u=1}^{n}
A_u^{(n+1),0}
\sum_{b \in \mathbb{F}} (b-r)^{(n+1)-u} Z(\mathbf{c},b,u) \\
& \quad + \sum_{u=n+2}^{k-2} A_u^{(n+1)} \sum_{b \ne r}
(b-r)^{(n+1)-u} Z(\mathbf{c},b,u).
\end{align*}}
Here $A_u^{(n+1),0} = \alpha_u(n+1,P), 1 \le u \le n+1,$ depend only
on $u$ and $n$.  By Proposition \ref{intersum}, the first triple sum
of the right hand side is equal to
\[
\sum_{u=1}^{n} \tfrac{\alpha_0(n+1,P)}{\lambda - \lambda_0(P)}
\sum_{s \in \mathbb{F}} (-1)^{u+1}\tbinom{n+1}{u}(s-r)^{(n+1)-u}
Z(\mathbf{c},s,u),
\]
which can be combined with the middle double sum of the right hand
side to bring the above identity to the following form:
\begin{align*}
(\lambda - &\lambda_{n+1}(P)) Z(\mathbf{c},r, n+1) \\
& = \sum_{u=1}^n A_u^{(n+1),1} \sum_{b \in \mathbb{F}}
(b-r)^{(n+1)-u} Z(\mathbf{c},b,u) + \sum_{u=n+2}^{k-2} A_u^{(n+1)}
\sum_{b \ne r} (b-r)^{(n+1)-u} Z(\mathbf{c},b,u).
\end{align*}
 Next we
replace $Z(\mathbf{c},b,1)$ above by (\ref{indhypo})$_1$ and use
Proposition \ref{intersum} to express $(\lambda - \lambda_{n+1}(P))$
times $Z(\mathbf{c},r, n+1)$ as a linear combination of $\sum_{b \in
\mathbb{F}} (b-r)^{n+1-u} Z(\mathbf{c},b,u)$ for $2 \le u \le n$ and
$\sum_{b \ne r} (b-r)^{n+1-u} Z(\mathbf{c},b,u)$ for $n+2 \le u \le
k-2$ with coefficients $A_u^{(n+1),2}$ depending only on $n$ and
$u$. Repeat this procedure. After $n-1$ iterations, we arrive at
\begin{align*}
(\lambda - \lambda_{n+1}(P))& Z(\mathbf{c},r, n+1) \\
&= A_n^{(n+1),n} \sum_{b \in \mathbb{F}} (b-r) Z(\mathbf{c},b,n) +
\sum_{u=n+2}^{k-2} A_u^{(n+1)} \sum_{b \ne r} (b-r)^{(n+1)-u}
Z(\mathbf{c},b,u).
\end{align*}
For the final calculation, use (\ref{indhypo})$_n$ to get
\begin{align*}
(&\lambda - \lambda_{n+1}(P)) Z(\mathbf{c},r, n+1) \\
&= \tfrac{A_n^{(n+1),n}}{\lambda - \lambda_n(P)} \sum_{b \in
\mathbb{F}} \sum_{u=n+1}^{k-2} A_u^{(n)} \sum_{s \ne b}
\tfrac{b-r}{(s-b)^{u-n}} Z(\mathbf{c},s,u) + \sum_{u=n+2}^{k-2}
A_u^{(n+1)}
\sum_{b \ne r} (b-r)^{(n+1)-u} Z(\mathbf{c},b,u) \\
&= \sum_{u=n+2}^{k-2} A_u^{(n+1)} \sum_{b \ne r} (b-r)^{(n+1)-u}
Z(\mathbf{c},b,u).
\end{align*}
Hence Lemma \ref{red} follows by induction.
\end{proof}

The techniques used to prove Lemma \ref{red} can be extended to
describe the eigenspaces of $T_\mathfrak{P}$. Let $\mathbf{c}$ be an
eigenfunction of $T_\mathfrak{P}$ with eigenvalue $\lambda_n(P)$.
The relations among $Z(\mathbf{c},r,j)$ for $r \in \mathbb F$ and $0
\le j \le k-2$ are distinguished by two cases, according to
$\lambda_j(P)$ equal to $\lambda_n(P)$ or not.

For those $l$ with $\lambda_l(P) \ne \lambda_n(P)$, the equation
(\ref{pflemma})$_l$ gives
\begin{align*}
Z(\mathbf{c},b,l) = \tfrac1{\lambda_n(P) - \lambda_l(P)} \Bigl[
&\sum_{u=0}^{l-1} \alpha_u(l,P) \sum_{s \in \mathbb{F}}
(s-b)^{l-u} Z(\mathbf{c},s,u) \\
&-\sum_{u=l+1}^{k-2} \beta_u(l, P) \sum_{s \ne b} (s-b)^{l-u}
Z(\mathbf{c},s,u) \Bigr]
\end{align*}
for all $b \in \mathbb{F}$. Further, for such $l$ we have $\sum_{b
\in \mathbb F} Z(\mathbf{c},b,l) = 0$ by (\ref{sumzero}). Let $l_0 <
l_1 < \dots < l_t$ be the distinct $l$'s such that $\lambda_l(P) \ne
\lambda_n(P)$. Then the same inductive procedure as in the proof of
Lemma \ref{red} yields, for each $l_v$, $0 \le v \le t$,
\begin{subequations} \label{nequ}
\begin{equation} \tag*{(\ref{nequ}$)_{l_v}$}
Z(\mathbf{c},b,l_v) = \sum_{\substack{0 \le u < l_v \\
\lambda_u(P) = \lambda_n(P)}} A_u^{(l_v)} \sum_{s \in \mathbb{F}}
(s-b)^{l_v-u} Z(\mathbf{c},s,u) + \sum_{u=l_v+1}^{k-2} A_u^{(l_v)}
\sum_{s \ne b} (s-b)^{l_v-u} Z(\mathbf{c},s,u)
\end{equation}
\end{subequations}
for some explicitly determined elements $A_u^{(l_v)}$ in
$\mathbb{F}(T)$ depending only on $u$ and $P$.

Let $i$ be an index such that $\lambda_i(P) = \lambda_n(P)$. The
Hecke action (\ref{redhegam}) gives rise to
\begin{equation} \label{equ}
\begin{aligned}
0 = & \sum_{u=0}^{i-1} \alpha_u(i, P) \sum_{b \in \mathbb{F}}
(b-r)^{i-u} Z(\mathbf{c},b,u)  +  [P^{k-2-i} - \lambda_i(P)] \sum_{b
\in
\mathbb{F}} Z(\mathbf{c},b,i) \\
&- \sum_{u=i+1}^{k-2} \beta_u(i, P) \sum_{b \ne r} (b-r)^{i-u}
Z(\mathbf{c},b,u).
\end{aligned}
\end{equation}
By successively substituting (\ref{nequ})$_{l_v}$ into (\ref{equ}),
starting with $v = 0$ and ending with $v = t$, and simplifying the
expression using Proposition \ref{intersum} at each step, we
eliminate all $Z(\mathbf{c},b, l)$'s in the equation (\ref{equ})
with $\lambda_l(P) \ne \lambda_n(P)$ and arrive at an identity of
the form
\begin{subequations} \label{finaleli}
\begin{equation} \tag*{(\ref{finaleli})$_{i,r}$}
0 = \sum_{\substack{0 \le u \le k-2 \\
\lambda_u(P) = \lambda_n(P)}} C_u(i,P) \sum_{b \ne r} (b-r)^{i-u}
Z(\mathbf{c},b,u)
\end{equation}
\end{subequations}
for some explicitly determined elements $C_u(i,P)$ in
$\mathbb{F}(T)$ depending only on $i, u$ and $P$. We have shown

\begin{thm} \label{mainga}
Suppose $q \ge k \ge 2$. Let $\mathfrak{P} = (P)$, where $P \in
\mathbb F[T]$ has degree one and $P(0)=1$. Then $\lambda_i(P)$, $0
\le i \le k-2$, with suitable multiplicities are the eigenvalues of
the Hecke operator $T_\mathfrak{P}$ on $S_{k,m}(\Gamma(T))$. For $0
\le n \le k-2$, set $A_n = \{ i: 0 \le i \le k-2 \; \text{and}\;
\lambda_i (P) = \lambda_n(P) \}$ and denote the integers in $[0,
~k-2]\smallsetminus A_n$ by $ l_0 < \cdots < l_t.$ Let $\mathbf{c}$
be an eigenfunction in $H_{k,m}(\Gamma(T))$ with eigenvalue
$\lambda_n(P)$. Then $\mathbf{c}$ is determined by
$Z(\mathbf{c},b,u)$ with $u \in A_n$ and $b \in \mathbb{F}$ subject
to the conditions $(\ref{finaleli})_{i,r}$ for $i \in A_n$ and $r
\in \mathbb{F}$. The remaining $Z(\mathbf{c},b,l)$'s are determined
recursively by $(\ref{nequ})_{l_v}$ from $v=t$ to $v=0$.
\end{thm}

\medskip

\section{Examples}

To illustrate Theorem \ref{mainga}, we compute the action of
$T_\mathfrak{P}$ on $H_{k,m}(\Gamma)$ for small weights $k=3, 4,
5$. None of these are diagonalizable with respect to the Hecke
operator. Let $\mathbf{c}$ be an eigenfunction.
\smallskip

\noindent\textbf{(i)} $q \ge k=3$. Here $\lambda_0(P) = \lambda_1(P)
= 1$. It follows from (\ref{redhegam}) that
\[
\sum_{b \in \mathbb{F}} Z(\mathbf{c},b,0) = \sum_{b \ne r}
\frac{Z(\mathbf{c},b,1)}{r-b}
\]
for all $r \in \mathbb{F}$. We shall solve this linear system. Fix
a generator $a$ of $\mathbb{F}^\times$ and arrange the elements of
$\mathbb{F}$ in the order $0, a, a^2, \dots, a^{q-1}$. Express the
above system in matrix form
\begin{equation} \label{k3sys}
\left(\begin{matrix}
  0 & -\tfrac{1}{a} & -\tfrac{1}{a^2}  & -\tfrac{1}{a^3} & \dots &
  -\tfrac{1}{a^{q-1}} \\
  \tfrac{1}{a} & 0 & \tfrac{1}{a-a^2}  & \tfrac{1}{a-a^3} & \dots &
  \tfrac{1}{a-a^{q-1}} \\
  \tfrac{1}{a^2} & \tfrac{1}{a^2-a}  & 0 &  \tfrac{1}{a^2-a^3} & \dots &
  \tfrac{1}{a^2-a^{q-1}} \\
  \vdots & \vdots & \vdots & \vdots & \ddots & \vdots \\
  \tfrac{1}{a^{q-1}} & \tfrac{1}{a^{q-1}-a}  & \tfrac{1}{a^{q-1}-a^2} &
  \tfrac{1}{a^{q-1}-a^3} & \dots & 0
\end{matrix}\right)\left(\begin{matrix}
  Z(\mathbf{c},0,1) \\ Z(\mathbf{c},a,1) \\ Z(\mathbf{c},a^2,1) \\ \vdots \\ Z(\mathbf{c},a^{q-1},1)
\end{matrix}\right) = \left(\begin{matrix}
  c \\ c \\ c \\ \vdots \\ c
\end{matrix}\right),
\end{equation}
where $c = \sum_{b \in \mathbb{F}} Z(\mathbf{c},b,0)$. We
determine the nullity of the coefficient matrix $M$. Write
\[
M = \left(\begin{matrix}
  1 & 0 & 0 & \dots & 0 \\ 0 & \tfrac{1}{a} & 0 & \dots & 0 \\ 0 & 0 & \tfrac{1}{a^2} & \dots & 0 \\
  \vdots & \vdots & \vdots & \ddots & \vdots \\ 0 & 0 & 0 & \dots &
  \tfrac{1}{a^{q-1}}
\end{matrix}\right)\left(\begin{matrix}
  0 & -\tfrac{1}{a} & -\tfrac{1}{a^2}  & -\tfrac{1}{a^3} & \dots &
  -\tfrac{1}{a^{q-1}} \\
  1 & 0 & \tfrac{1}{1-a}  & \tfrac{1}{1-a^2} & \dots &
  \tfrac{1}{1-a^{q-2}} \\
  1 & \tfrac{1}{1-a^{q-2}}  & 0 &  \tfrac{1}{1-a} & \dots &
  \tfrac{1}{1-a^{q-3}} \\
  \vdots & \vdots & \vdots & \vdots & \ddots & \vdots \\
  1 & \tfrac{1}{1-a}  & \tfrac{1}{1-a^2} &
  \tfrac{1}{1-a^3} & \dots & 0
\end{matrix}\right).
\]
Call the second matrix on the right hand side $C$. Note that
$\Null(M) = \Null(C)$. Consider the submatrix obtained from $C$ by
deleting the first row and the first column

\[
C' = \left(\begin{matrix}
  0 & \tfrac{1}{1-a}  & \tfrac{1}{1-a^2} & \dots &
  \tfrac{1}{1-a^{q-2}} \\
  \tfrac{1}{1-a^{q-2}}  & 0 &  \tfrac{1}{1-a} & \dots &
  \tfrac{1}{1-a^{q-3}} \\
  \vdots & \vdots & \vdots & \ddots & \vdots \\
  \tfrac{1}{1-a}  & \tfrac{1}{1-a^2} & \tfrac{1}{1-a^3} & \dots & 0
\end{matrix}\right),
\]
which is a $(q-1)\times(q-1)$ circulant matrix. Then $\mathbf{v}'_j
= \left(\begin{smallmatrix} 1 \\ a^j \\ a^{2j} \\ \vdots \\
a^{(q-2)j}
\end{smallmatrix}\right)$, $j = 1, 2, \dots, q-1$, are $q-1$
linearly  independent eigenvectors of $C'$ with eigenvalue
\[
\frac{a^j}{1-a} + \frac{a^{2j}}{1-a^2} + \cdots +
\frac{a^{(q-2)j}}{1-a^{q-2}} = j,
\]
as a consequence of the following lemma.

\begin{lem} \label{binomsum}
For $j=1, 2, \dots, q-1$ and $l \ge 1$, we have
\[
\sum_{n=1}^{q-2} \frac{a^{jn}}{(1-a^n)^l} = (-1)^{l-1}\tbinom{j}{l}.
\]
\end{lem}
\begin{proof}
We shall prove this lemma by induction on $l$. For $l=1$, we compute
\[
\sum_{n=1}^{q-2} \frac{a^{jn}}{1-a^n} = \sum_{n=1}^{q-2}
\frac{a^{jn}-1+1}{1-a^n} = \sum_{n=1}^{q-2}\Bigl[- (1 + a^n + \cdots
+ a^{(j-1)n}) + \frac{1}{1-a^n} \Bigr].
\]
Since $a$ has order $q-1$ and $1 \le j \le q-1$, $\sum_{n=1}^{q-2}
a^{in} = -a^{i(q-1)} = -1$ for all $i = 1, \dots, j-1$. As
$\sum_{n=1}^{q-2} \frac{1}{1-a^n} = -1$, the above sum is equal to
\[
\sum_{n=1}^{q-2} \frac{a^{jn}}{1-a^n} = -(q-2)+(j-1)-1 =j.
\]
Next, we assume that $\sum_{n=1}^{q-2} \frac{a^{jn}}{(1-a^n)^l} =
(-1)^{l-1}\tbinom{j}{l}$ for all $j = 1, \dots, q-1$. Then
\begin{align*}
\sum_{n=1}^{q-2} \frac{a^{jn}}{(1-a^n)^{l+1}} &= \sum_{n=1}^{q-2}
\frac{a^{jn}-1+1}{(1-a^n)^{l+1}} = \sum_{n=1}^{q-2} \Bigl[ -
\frac{1+a+\dots+a^{(j-1)n}}{(1-a^n)^l} + \frac{1}{(1-a^n)^{l+1}}
\Bigr]\\
&= -[-1 + (-1)^{l-1}\tbinom{1}{l} + \cdots +
(-1)^{l-1}\tbinom{j-1}{l}] - 1 = (-1)^l\tbinom{j}{l+1}
\end{align*}
by the Pascal's triangle identity $\sum_{i=1}^m \tbinom{i}{l} =
\tbinom{m+1}{l+1}$. The lemma follows by induction.
\end{proof}

Back to the matrix $C$. The vectors $\mathbf{v}_0 =
\left(\begin{smallmatrix} 1 \\ 1 \\ 1 \\ \vdots \\ 1
\end{smallmatrix}\right)$ and $\mathbf{v}_j
= \left(\begin{smallmatrix} 0 \\ 1 \\ a^j \\ \vdots \\
a^{(q-2)j}
\end{smallmatrix}\right)$, $j = 1, \dots, q-1$, are $q$ linearly
independent eigenvectors of $C$ with the eigenvalues $0$ and $j$,
respectively. Since our field has characteristic $p$, this shows
that the $\Nul(C) = q/p$. If $c = \sum_{b \in \mathbb{F}}
Z(\mathbf{c},b,0) = 0$, then we obtain $(q-1)+{q}/{p}$ linearly
independent eigenvectors for $T_\mathfrak{P}$. When $c \ne 0$, note
that $\mathbf{v}
= \left(\begin{smallmatrix} 0 \\ ca \\ ca^2 \\ \vdots \\
ca^{q-1}
\end{smallmatrix}\right)$ is a solution of (\ref{k3sys}).
Together with the homogeneous ones, we have $q+{q}/{p}$ linearly
independent eigenvectors for $T_\mathfrak{P}$, all with eigenvalue
$1$. Since $1$ is the only eigenvalue of $T_\mathfrak{P}$, its
total multiplicity $2q$, thus $T_\mathfrak{P}$ is not
diagonalizable. We record this result in

\begin{prop}
Suppose $\mathbb F$ has cardinality $q \ge 3$ and characteristic
$p$. For a maximal degree one ideal $\mathfrak{P}\ne(T)$, $1$ is
the only eigenvalue of the Hecke operator $T_{\mathfrak{P}}$ on
$S_{3,m}(\Gamma(T))$. The eigenspace of $T_\mathfrak{P}$ has
dimension $q+q/p$, while the space $S_{3,m}(\Gamma(T))$ has
dimension $2q$. Consequently, $T_{\mathfrak{P}}$ is not
diagonalizable on $S_{3,m}(\Gamma(T))$.
\end{prop}

As the dimension of the $1$-eigenspace of $T_{\mathfrak{P}}$ on
$S_{3,m}^2(\Gamma(T))$ is $q-1$, which is $\dim_C
S_{3,m}^2(\Gamma(T))$, so $T_{\mathfrak{P}}$ is diagonalizable on
$S_{3,m}^2(\Gamma(T))$.

\medskip

\noindent\textbf{(ii)} $q \ge k=4$. In this case $\lambda_0(P) =
\lambda_2(P) = 1$ and $\lambda_1(P) = -P+2$. A similar computation
yields

\begin{prop}
Suppose $\mathbb F$ has cardinality $q \ge 4$ and characteristic
$p$. For a maximal degree one ideal $\mathfrak{P}\ne(T)$, $1$ and
$2-P$ are the two distinct eigenvalues of the Hecke operator
$T_{\mathfrak{P}}$ on $S_{4,m}(\Gamma(T))$. The $1$-eigenspace has
dimension $q+2q/p$ if $p>2$ and dimension $q+q/p$ if $p=2$. The
$(2-P)$-eigenspace has dimension $q$. Moreover, $T_{\mathfrak{P}}$
is not diagonalizable on $S_{4,m}(\Gamma(T))$.
\end{prop}

One checks that $T_{\mathfrak{P}}$ on $S_{4,m}^2(\Gamma(T))$ is
diagonalizable since $\dim_C ~S_{4,m}^2(\Gamma(T))= 2q-1$, the
$1$-eigenspace
 is $(q-1)$-dimensional, and the $(2-P)$-eigenspace has dimension
 $q$.
\medskip

\noindent\textbf{(iii)} $q \ge k=5$. In this case $ \lambda_0(P) =
\lambda_3(P) = 1$ and $\lambda_1(P) = \lambda_2(P) = -2P+3$. First
we assume $p > 2$ so that $1 \ne -2P+3$. To determine the
$1$-eigenspace, consider the equations from (\ref{redhegam}) with $j
= 0,1,2$: {\small\begin{equation} \label{k5j0} 0 = (P^2+P+1) \sum_{b
\in \mathbb{F}} Z(\mathbf{c},b,0) + 3 \sum_{b \ne r}
\frac{Z(\mathbf{c},b,1)}{b-r} - 3(P-1) \sum_{b \ne r}
\frac{Z(\mathbf{c},b,2)}{(b-r)^2} + (P-1)^2 \sum_{b \ne r}
\frac{Z(\mathbf{c},b,3)}{(b-r)^3},
\end{equation}}
\begin{align*}
2 Z(\mathbf{c},r,1) = & \:(P+1) \sum_{b \in \mathbb{F}} (b-r)
Z(\mathbf{c},b,0) +
(P+3) \sum_{b \in \mathbb{F}} Z(\mathbf{c},b,1) \\
&\: -(P-3) \sum_{b \ne r} \frac{Z(\mathbf{c},b,2)}{b-r} - (P-1)
\sum_{b \ne r} \frac{Z(\mathbf{c},b,3)}{(b-r)^2},
\end{align*}
and
\[
2Z(\mathbf{c},r,2) = \sum_{b \in \mathbb{F}} (b-r)^2
Z(\mathbf{c},b,0) + 3 \sum_{b \in \mathbb{F}} (b-r)
Z(\mathbf{c},b,1) +3 \sum_{b \in \mathbb{F}} Z(\mathbf{c},b,2) +
\sum_{b \ne r} \frac{Z(\mathbf{c},b,3)}{b-r}
\]
for all $r \in \mathbb{F}$. Summing the second and third equations
over all $r$, we get $\sum_{r \in \mathbb{F}} Z(\mathbf{c},r,1) =
0$ and $\sum_{r \in \mathbb{F}} Z(\mathbf{c},r,2) = 0$, which lead
to the following simplifications of the second and third
equations:
\begin{equation} \label{k5j1}
Z(\mathbf{c},r,1) = \frac{P+1}{2} \sum_{b \in \mathbb{F}} (b-r)
Z(\mathbf{c},b,0) - \frac{P-3}{2} \sum_{b \ne r}
\frac{Z(\mathbf{c},b,2)}{b-r} - \frac{P-1}{2} \sum_{b \ne r}
\frac{Z(\mathbf{c},b,3)}{(b-r)^2},
\end{equation}
and
\begin{equation} \label{k5j2'}
Z(\mathbf{c},r,2) = \frac{1}{2} \sum_{b \in \mathbb{F}} (b-r)^2
Z(\mathbf{c},b,0) + \frac{3}{2} \sum_{b \in \mathbb{F}} (b-r)
Z(\mathbf{c},b,1) + \frac{1}{2} \sum_{b \ne r}
\frac{Z(\mathbf{c},b,3)}{b-r}
\end{equation}
for all $r \in \mathbb{F}$. Plugging (\ref{k5j1}) into
(\ref{k5j2'}) and using Proposition \ref{intersum} to simplify, we
get
\begin{equation} \label{k5j2}
Z(\mathbf{c},r,2) = \frac{1}{2} \sum_{b \in \mathbb{F}} (b-r)^2
Z(\mathbf{c},b,0) + \frac{1}{2} \sum_{b \ne r}
\frac{Z(\mathbf{c},b,3)}{b-r}
\end{equation}
for all $r \in \mathbb{F}$. Substituting (\ref{k5j1}) and
(\ref{k5j2}) into (\ref{k5j0}) and simplifying the result using
Proposition \ref{intersum}, we obtain
\begin{equation} \label{k5l1}
\sum_{b \in \mathbb{F}} Z(\mathbf{c},b,0) =  \sum_{b \ne r}
\frac{Z(\mathbf{c},b,3)}{(r-b)^3}
\end{equation}
for all $r \in \mathbb{F}$. To solve the above linear system, we
employ the same method as in case (i), that is, computing the
nullity of
\[
C = \left(\begin{matrix}
  0 & -\tfrac{1}{a^3} & -\tfrac{1}{a^6}  & -\tfrac{1}{a^9} & \dots &
  -\tfrac{1}{a^{3(q-1)}} \\
  1 & 0 & \tfrac{1}{(1-a)^3}  & \tfrac{1}{(1-a^2)^3} & \dots &
  \tfrac{1}{(1-a^{q-2})^3} \\
  1 & \tfrac{1}{(1-a^{q-2})^3}  & 0 &  \tfrac{1}{(1-a)^3} & \dots &
  \tfrac{1}{(1-a^{q-3})^3} \\
  \vdots & \vdots & \vdots & \vdots & \ddots & \vdots \\
  1 & \tfrac{1}{(1-a)^3}  & \tfrac{1}{(1-a^2)^3} &
  \tfrac{1}{(1-a^3)^3} & \dots & 0
\end{matrix}\right).
\]
By Lemma \ref{binomsum}, the vectors $\mathbf{v}_0 =
\left(\begin{smallmatrix} 1 \\ 1 \\ 1 \\ \vdots \\ 1
\end{smallmatrix}\right)$ and $\mathbf{v}_j
= \left(\begin{smallmatrix} 0 \\ 1 \\ a^j \\ \vdots \\
a^{(q-2)j} \end{smallmatrix}\right)$, $j = 1, \dots, q-1$, are $q$
linearly independent eigenvectors of $C$ with the eigenvalues $0$
and $\binom{j}{3}$, respectively. Therefore the nullity of $C$ is
$3q/p$ when $p>3$ and is $q/p$ when $p=3$, which yields the number
of linearly independent eigenvectors if $c := \sum_{b \in
\mathbb{F}} Z(\mathbf{c},b,0) = 0$. When $c \ne 0$, we note that
$\mathbf{v}
= \left(\begin{smallmatrix} 0 \\ ca^3 \\ ca^6 \\ \vdots \\
ca^{3(q-1)} \end{smallmatrix}\right)$ is a solution of
(\ref{k5l1}). Together with the homogeneous ones, we see that the
$1$-eigenspace of $T_\mathfrak{P}$ has dimension $q+3q/p$ if $p
> 3$ and $q+q/p$ if $p=3$.

Next we determine the eigenvectors with eigenvalue $-2P+3$. Such
eigenvectors are double cusp forms by Proposition \ref{eigenone},
 so $Z(\mathbf{c},r,3) =0$ for all $r \in \mathbb{F}$. Thus the equations
 from
(\ref{redhegam}) with $j = 0, 1, 2$  can be simplified as
\begin{align*}
&Z(\mathbf{c},r,0) = - \frac{3}{2} \sum_{b \ne r}
\frac{Z(\mathbf{c},b,1)}{b-r} +
\frac{3(P-1)}{2} \sum_{b \ne r} \frac{Z(\mathbf{c},b,2)}{(b-r)^2}, \\
&0 = (P+1)\sum_{b \in \mathbb{F}} (b-r) Z(\mathbf{c},b,0) + (P+3)
\sum_{\mathbb{F}} Z(\mathbf{c},b,1) - (P-3) \sum_{b \ne r}
\frac{Z(\mathbf{c},b,2)}{b-r},
\end{align*}
and
\[
0 = \sum_{b \in \mathbb{F}} (b-r)^2 Z(\mathbf{c},b,0) + 3 \sum_{b
\in \mathbb{F}} (b-r) Z(\mathbf{c},b,1) + 3 \sum_{b \in \mathbb{F}}
Z(\mathbf{c},b,2)
\]
for all $r \in \mathbb{F}$. Substituting the first relation into
the second and the third, and simplifying the resulting
expressions by using Proposition \ref{intersum}, we arrive at
\begin{equation}\label{k5lp1}
\sum_{b \in \mathbb{F}} Z(\mathbf{c},b, 1) = 2 \sum_{b \ne r}
\frac{Z(\mathbf{c},b,2)}{r-b}
\end{equation}
and
\begin{equation}\label{k5lp}
\sum_{b \in \mathbb{F}} Z(\mathbf{c},b,2) = 0
\end{equation}
for all $r \in \mathbb{F}$. Write $c = \sum_{b \in \mathbb{F}}
Z(\mathbf{c},b, 1)$. Solve the system (\ref{k5lp1}) using the same
method as (\ref{k3sys}). When $c = 0$, we get homogeneous solutions
$\left(\begin{smallmatrix} Z(\mathbf{c},0,2) \\ Z(\mathbf{c},a,2) \\ Z(\mathbf{c},a^2,2) \\ \vdots \\
Z(\mathbf{c},a^{q-1},2) \end{smallmatrix}\right) = \left(\begin{smallmatrix} 1 \\
1 \\ 1 \\ \vdots \\ 1 \end{smallmatrix}\right) \: \text{or} \:
\left(\begin{smallmatrix} 0 \\ 1 \\ a^j \\ \vdots \\ a^{(q-2)j}
\end{smallmatrix}\right)$, $j=1, \dots, q-1$; when $c \ne 0$, we
get a nonhomogeneous solution
$\left(\begin{smallmatrix} Z(\mathbf{c},0,2) \\ Z(\mathbf{c},a,2) \\ Z(\mathbf{c},a^2,2) \\ \vdots \\
Z(\mathbf{c},a^{q-1},2) \end{smallmatrix}\right) = \frac{1}{2}\left(\begin{smallmatrix} 0 \\
ca \\ ca^2 \\ \vdots \\ ca^{q-1} \end{smallmatrix}\right)$. Note
that all solutions satisfy the equation (\ref{k5lp}). Thus the
$(-2P+3)$-eigenspace of $T_\mathfrak{P}$ has dimension $q+q/p$.
Combined with the dimension of $1$-eigenspace, we conclude that
$T_\mathfrak{P}$ is not diagonalizable on $S_{5,m}(\Gamma(T))$ since
the space has dimension $4q$. We summarize the above discussion in

\begin{prop}
Suppose $\mathbb F$ has cardinality $q \ge 4$ and characteristic
$p>2$. For a maximal degree one ideal $\mathfrak{P}\ne(T)$, $1$
and $-2P+3$ are the two distinct eigenvalues of the Hecke operator
$T_{\mathfrak{P}}$ on $S_{5,m}(\Gamma(T))$. The $1$-eigenspace has
dimension $q+3q/p$ if $p>3$ and dimension $q+q/p$ if $p=3$. The
$(-2P+3)$-eigenspace has dimension $q+q/p$. Further,
$T_{\mathfrak{P}}$ is not diagonalizable on $S_{5,m}(\Gamma(T))$.
\end{prop}

Now we turn to the case when $p=2$. In this case, we have only one
eigenvalue, namely, 1.  Then (\ref{redhegam}) for $j=0,1,2$ become
{\small\begin{align*} 0 &=(P^2+P+1) \sum_{b \in \mathbb{F}}
Z(\mathbf{c},b,0)+ \sum_{b \ne r} \frac{Z(\mathbf{c},b,1)}{b-r} +
(P-1) \sum_{b \ne r} \frac{Z(\mathbf{c},b,2)}{(b-r)^2} +
(P-1)^2 \sum_{b \ne r} \frac{Z(\mathbf{c},b,3)}{(b-r)^3}, \\
0 &= \sum_{b \in \mathbb{F}} (b-r) Z(\mathbf{c},b,0) + \sum_{b \in
\mathbb{F}} Z(\mathbf{c},b,1) + \sum_{b \ne r}
\frac{Z(\mathbf{c},b,2)}{b-r} + \sum_{b \ne r}
\frac{Z(\mathbf{c},b,3)}{(b-r)^2},
\end{align*}}
and
\[
0 = \sum_{b \in \mathbb{F}} (b-r)^2 Z(\mathbf{c},b,0) + \sum_{b \in
\mathbb{F}} (b-r) Z(\mathbf{c},b,1) + \sum_{b \in \mathbb{F}}
Z(\mathbf{c},b,2) + \sum_{b \ne r} \frac{Z(\mathbf{c},b,3)}{(b-r)}
\]
for all $r \in \mathbb{F}$. Observe that we can represent the
above system as a homogeneous matrix equation $M \mathbf{x} =
\mathbf{0}$, where $M$ is a $3q \times 4q$ matrix. Moreover, it is
clear that rank$\,M > 1$. Thus the eigenspace is has dimension
less than $4q$, so that the Hecke operator $T_\mathfrak{P}$ is not
diagonalizable. Therefore we have shown

\begin{prop}
Suppose $\mathbb F$ has cardinality $q \ge 4$ and characteristic
$p=2$. For a maximal degree one ideal $\mathfrak{P}\ne(T)$, $1$
 is the only eigenvalue of the Hecke operator
$T_{\mathfrak{P}}$ on $S_{5,m}(\Gamma(T))$. The eigenspace of
$T_\mathfrak{P}$ has dimension less than $4q$, the dimension of
$S_{5,m}(\Gamma(T))$. Hence $T_{\mathfrak{P}}$ is not
diagonalizable on $S_{5,m}(\Gamma(T))$.
\end{prop}

As for the action of $T_{\mathfrak{P}}$ on $S_{5,m}^2(\Gamma(T))$,
by the same computation as before, we see that for $q$ odd, the
$1$-eigenspace is $(q-1)$-dimensional and the $(3-2P)$-eigenspace
has dimension $q + q/p$ so that the total dimension is less than
$3q-1$, the dimension of $S_{5,m}^2(\Gamma(T))$; for $q$ even, the
matrix $M$ is $3q \times 3q$ with rank at least two, thus the
eigenspace  is at most $(3q - 2)$-dimensional. Hence in both
cases, $T_{\mathfrak{P}}$ on $S_{5,m}^2(\Gamma(T))$ is not
diagonalizable.

\begin{remark}
For Drinfeld cusp forms, what happens in case (iii) is
representative of the general weights.  For example, when the weight
$k=6$, we have three distinct eigenvalues $1, 4-3P$ and $6-6P+P^2$
if $p \ne 3$ and two distinct eigenvalues $1$ and $P^2$ if $p=3$.
The computations for $Z(\mathbf{c},b,u)$ are similar.
\end{remark}

\end{document}